\documentclass[12pt]{article}

\usepackage[left=2.5cm, right=2.5cm, top=2cm, bottom=2cm]{geometry}

\usepackage{bbm}
\usepackage{amsmath, amsthm, amssymb}
\usepackage{amsfonts}
\usepackage{microtype}
\usepackage{pdfrender}
\usepackage[ansinew]{inputenc}
\usepackage[dvips]{epsfig}
\usepackage{graphicx}
\usepackage[english]{babel}

\usepackage{cite}
\usepackage{graphicx}
\usepackage{amscd}
\usepackage{xcolor}
\usepackage{bm}
\usepackage{enumerate}
\everymath{\displaystyle}
\usepackage{verbatim}
\usepackage{hyperref}
\usepackage{amstext}
\usepackage{latexsym}
\let\oldsqrt\sqrt
\def\sqrt{\mathpalette\DHLhksqrt}
\def\DHLhksqrt#1#2{%
\setbox0=\hbox{$#1\oldsqrt{#2\,}$}\dimen0=\ht0
\advance\dimen0-0.2\ht0
\setbox2=\hbox{\vrule height\ht0 depth -\dimen0}%
{\box0\lower0.4pt\box2}}

\allowdisplaybreaks

\newcommand{\R}{\mathbb{R}} 
\newcommand{\N}{\mathbb{N}} 

\newcommand{\sign}{\textnormal{sign}} 
\newcommand{\supp}{\textnormal{supp}} 

\renewcommand{\div}{\textnormal{div}}

\newcommand{\cC}{{\mathcal C}}

\newcommand{\cE}{{\mathcal E}}

\newcommand{\cN}{{\mathcal N}}

\newcommand{\cQ}{{\mathcal Q}}

\newcommand{\eps}{\varepsilon}

\theoremstyle{definition}
\newtheorem{defi}{Definition}[section]
\newtheorem{remark}[defi]{Remark}

\theoremstyle{plain} 
\newtheorem{thm}[defi]{Theorem}
\newtheorem{prop}[defi]{Proposition}
\newtheorem{lemma}[defi]{Lemma}

\theoremstyle{definition}

\numberwithin{equation}{section}

 \title{ Eigenvalues   of nonlinear $(p,q)$-fractional Laplace operators  under nonlocal Neumann  conditions }

\author{Pierre Aime Feulefack $^a$\footnote{Email: \href{mailto:feulefac@sas.upenn.edu}{feulefac@sas.upenn.edu}} ~ and ~ Emmanuel Wend-Benedo Zongo \vspace{.3cm} $^b$\footnote{Corresponding author. Email: \href{mailto:emmanuel.zongo@math.univ-toulouse.fr}{emmanuel.zongo@math.univ-toulouse.fr} } 
\\
\footnotesize{$^a$ Department of Mathematics, University of Pennsylvania,}\\
\footnotesize{David Rittenhouse Laboratory, 209 S 33rd St, Philadelphia, PA 19104, US\vspace{.1cm}}
\\
\footnotesize{$^b$ Institut de math\'ematiques de Toulouse, Universit\'e Paul Sabatier, } \\
\footnotesize{  118, route de Narbonne F-31062 Toulouse Cedex 9, France }
}
 
\date{}

\begin{document}
\pdfrender{StrokeColor=black,TextRenderingMode=2,LineWidth=0.2pt}
\maketitle
\vspace{-.5cm}
\begin{abstract}
In this paper,  we investigate on a bounded open set of $\mathbb{R}^N$ with smooth boundary, an eigenvalue problem involving the sum of nonlocal operators $ (-\Delta)_p^{s_1}+ (-\Delta)_q^{s_2}$ with $s_1,s_2\in (0,1)$, $p,q\in (1,\infty)$ and subject to the corresponding homogeneous nonlocal $(p,q)$-Neumann boundary condition. 
 A careful analysis of the considered problem  leads us to a complete description of the set of eigenvalues as being the precise interval  $\{0\}\cup(\lambda_{1}(s_2,q),\infty)$, where $\lambda_{1}(s_2,q)$ is the first nonzero eigenvalue of the homogeneous fractional $q$-Laplacian under nonlocal $q$-Neumann boundary condition. Furthermore, we establish that every eigenfunctions   is globally bounded.
\end{abstract}

{\footnotesize
\begin{center}
\end{center}
\vspace{.4cm}
}

\tableofcontents

\section{Introduction and main results}
 In this article, we are interested in the study of a nonlinear nonlocal eigenvalue problem driven by the fractional $(p,q)$-Laplace operators subject to  homogeneous nonlocal $(p,q)$-Neumann boundary condition,
\begin{equation}\label{Eq1}
	\begin{split}
	\quad\left\{\begin{aligned}
		 (-\Delta)_p^{s_1}u+(-\Delta)_q^{s_2}u&= \lambda |u|^{q-2}u && \text{ in \quad $\Omega$}\\
		\cN_{p,s_1} u +\cN_{q,s_2} u&=0    && \text{ on }\quad  \mathbb{R}^N\setminus \overline{\Omega},
	\end{aligned}\right.
	\end{split}
	\end{equation}
where $\Omega\subset \R^N$ ($N\ge 2$) is an open bounded set  with  Lipschitz boundary $\partial\Omega$, and $p,q\in (1,\infty)$. 
The operator $(-\Delta)_p^{s}$, $s\in (0,1)$ stands for the fractional $p$-Laplacian and it is defined  up to a normalization constant and  for compactly supported smooth function $u:\R^N\to \R$ by
\begin{equation}
\begin{split}
 (-\Delta)_r^{s} u(x)&=\lim_{\varepsilon\to 0^+}\int_{\R^N\setminus B_{\varepsilon}(x)}\frac{|u(x)-u(y)|^{r-2}(u(x)-u(y))}{|x-y|^{N+rs}}\ dy.
 \end{split}
\end{equation}
This covers the usual definition of the fractional Laplacian $(-\Delta)^s$ when $r=2$ and the classical $r$-Laplacian $\Delta_r$ when $s=1$ with smooth $u$. The nonlocal $r$-Neumann condition   $\cN_{p,s} $, also called nonlocal normal $r$-derivative and describing the natural Neumann boundary condition in presence of the fractional $r$-Laplacian is defined by
\begin{equation}\label{Normalderivative}
\cN_{r,s} u(x)=C_{N,s}\int_{\Omega} \frac{
|u(x)-u(y)|^{r-2}(u(x)-u(y))}{|x-y|^{N+rs}}\ dy\qquad (x\in \R^N\setminus\overline\Omega).
\end{equation}
To a great extent, the study of nonlocal equations is motivated by real world applications. This is due to their ability to describe certain phenomena with large scale behavior with a possible better efficiency. In physics, they model phenomena like anomalous diffusion and long-range interactions in materials. In biology, they describe processes such as collective behavior in animal swarms and the spread of diseases, while in finance, they are used to capture market dynamics influenced by global factors. The operator $(-\Delta)_p^{s_1}+(-\Delta)_q^{s_2}$ called here fractional $(p,q)$-Laplacian, has been extensively investigated in recent years, under  Dirichlet boundary conditions.   We refer to \cite{BS22,FMS17,BM18,ZF24}, see also the references therein, for the Dirichlet case.
It is well known (see for example \cite{DS15})
,  that for smooth functions $u:\R^N\to \R$ with compact support,
$$\text{$(-\Delta)_r^{s}u\to -\Delta_ru$\quad  and \quad $\cN_{r,s} u\to|\nabla u|^{r-2}\partial_{\nu}u$\quad   as \quad$ s\to 1^-$.}$$
Therefore, in  the limiting case $s_1=s_2=1$, the fractional $(p,q)$-Laplacian reduces to the classical well-known $(p,q)$-Laplace operators $-\Delta_p-\Delta_q$, where for $\theta\in (1,\infty)$, formally, the operator $\Delta_{\theta}u :=\div(|\nabla u|^{\theta-2}\nabla u)$  stands for the classical  $\theta$-Laplace operator. 

The study of problems involving the $(p,q)$-Laplace operators  have received by far a lot of attention by many authors. This is due to their wide applications in physics and related sciences such as  biophysics, nonlinear optics, fluid mechanics, quantum and plasma physics, chemical reaction design \cite{BM23,BM2023,ZR24}.
The purely local counterpart of problem \eqref{Eq1} is when asymptotically, $s_1=s_2=1$ (see  \cite{BM19,MM16}). In that particular case, we have the following Neumann eigenvalue problem  involving the classical $(p,q)$-Laplace operators:
\begin{equation}\label{pq-Laplace}
	\begin{split}
	\quad\left\{\begin{aligned}
		 -\Delta_pu-\Delta_qu&= \lambda |u|^{q-2}u && \text{ in \quad $\Omega$}\\
		|\nabla u|^{p-2}\partial_{\nu}u&+|\nabla u|^{q-2}\partial_{\nu}u =0    && \text{ on }\quad  \partial\Omega,\\
	\end{aligned}\right.
	\end{split}
	\end{equation}
	where $\nu$ is the unit outward normal to $\partial\Omega $. Note that problem \eqref{Eq1} generalize in a suitably sense problem \eqref{pq-Laplace}, see also problem \eqref{Eq-Particular} for mixed local-nonlocal equations. Note that  problem \eqref{pq-Laplace} is also known in the literature as two phase problems and arises in other fields of sciences and Engineering, including, chemical reaction-diffusion system and nonlinear elasticity theory. We mention that the literature related to problems governed by the $(p, q)$-Laplace operators is vast,  have attracted considerable interest and  in fact, it is daily increasing.

Eigenvalue problems for  the $(p,q)$- Laplace operators like \eqref{pq-Laplace} have been considered by many authors, see for example  papers \cite{MM16} and also \cite{FMS15, BM23, BM2023, ZR24}  and the references in there, where, the authors  analysed the spectrum  of the operator  $-\Delta_q-\Delta_p$. In fact, it  has been shown that the set of eigenvalues of $\eqref{pq-Laplace}$ is precisely described by the interval $(\lambda_1(q),+\infty)$ plus an isolated point $\lambda=0$, where $\lambda_1(q)$ is the first positive eigenvalue of the $q$-Laplacian $-\Delta_q$ in $\Omega$, under the Neumann boundary condition.

The main goal of this paper is to  investigate the spectrum of problem \eqref{Eq1}, which  complements  the results of \cite{MM16}, \cite{FMS15} and \cite{BM023, BM23} from local to a nonlocal comprehensive analysis,  under nonlocal $(p,q)$-Neumann condition. In addition,  we establish that any eigenvalue of \eqref{Eq1} is globally bounded.  To our best knowledge, there is no literature related to   discussion for nonlinear eigenvalue problem  like \eqref{Eq1}, involving the fractional  $(p,q)$-Laplace  operators  $(-\Delta)_p^{s_1}+(-\Delta)_q^{s_2}$, under nonlocal $(p,q)$-Neumann condition.

The solution $u$ of problem \eqref{Eq1} is sought in the Sobolev  space $X_{\sigma,\theta}$ (see Section \ref{Preli} for definition), where the parameters $\sigma\in (0,1)$ and $\theta\in (1,\infty)$ are defined by $\sigma=\max\{s_1,s_2\}$ and $\theta =\max\{p,q\}$,  so the $(p,q)$-Neumann boundary condition exists in a  trace sense, and the problem is satisfied in the distribution sense (see \cite{DRV17, ML19}). Note that we only assume $p\neq q$ and allow various  possibilities for the parameters $s_1$ and $s_2$, see the case examples $(P_1)-(P_4)$ after Theorem \ref{Main-1} below.
Using the integration by parts formula for nonlocal Neumann problem, we have the following definition. 
\begin{defi} A number $\lambda\in\R$ is an eigenvalue of problem \eqref{Eq1} with $(p,q)$-Neumann boundary condition, if there exists $u_{\lambda}\in  X_{\sigma,\theta}\setminus\{0\}$ such
that
\begin{equation}\label{weakform}
\begin{array}{ll}
    &\iint_{\cQ}\frac{|u_{\lambda}(x)-u_{\lambda}(y)|^{p-2}(u_{\lambda}(x)-u_{\lambda}(y))(v(x)-v(y))}{|x-y|^{N+ps_1}}\ dxdy \vspace{10pt}\\
    &+\iint_{\cQ}\frac{|u_{\lambda}(x)-u_{\lambda}(y)|^{q-2}(u_{\lambda}(x)-u_{\lambda}(y))(v(x)-v(y))}{|x-y|^{N+qs_2}}\ dxdy \vspace{10pt}=\lambda\int_{\Omega}|u_{\lambda}|^{q-2}u_{\lambda}v\ dx,
    \end{array}
\end{equation}
for all $v\in X_{\sigma,\theta}$. Here and in the following,  we denote by $\Omega^c:=\R^N\setminus\Omega$ and 
$$\text{$\cQ:= (\R^N\times\R^N)\setminus(\Omega^c)^2 = (\Omega\times\Omega)\cup(\Omega\times\Omega^c)\cup(\Omega^c\times\Omega)$.}$$
\end{defi}
We shall call the function $u_{\lambda}\neq 0$ an eigenfunction of problem \eqref{Eq1} associated to $\lambda$, and $(\lambda,u_\lambda)$ is called eigenpairs of problem \eqref{Eq1}.
Note that taking $v=u_{\lambda}$ in \eqref{weakform}, we infer clearly that no negative  $\lambda$ can be an eigenvalue of problem \eqref{Eq1}. In addition, we observe that  $\lambda=0$ is an eigenvalue of problem \eqref{Eq1} with corresponding eigenfunctions being the nontrivial constants functions. So, we need only to investigate the case $\lambda>0$. Moreover, observe that if $\lambda>0$ is an eigenvalue of problem \eqref{Eq1}, then choosing  $v\equiv 1$ in \eqref{weakform}, we find that
$$\int_{\Omega}|u_{\lambda}|^{q-2}u_{\lambda}~dx=0.$$
Hence, the eigenfunctions corresponding to positive eigenvalues $\lambda$ of problem \eqref{Eq1} necessarily belong to the nonempty, symmetric, weakly closed cone
$$\mathcal{C}:=\{v\in X_{\sigma,\theta}~:~~\int_{\Omega}|v|^{q-2}v~dx=0\}.$$
Not that the set $\mathcal{C}\backslash\{0\}\subset X_{\sigma,\theta}$ is not empty. Indeed, let $u_1$ and $u_2$ be two nonnegative test functions having supports in two disjoint ball included in $\Omega$ such that $$\int_\Omega u^{q-1}_1dx=\int_\Omega u^{q-1}_2dx.$$ More precisely, let $x_1,x_2\in\Omega$ such that $x_1\neq x_2.$ Then there exists en $\eps>0$ small enough such that the balls $B_\eps(x_1)$ and  $B_\eps(x_2)$ are included in $\Omega$ and $B_\eps(x_1)\cap B_\eps(x_2)=\emptyset.$ Define the functions $u_j,$ $j=1,2$ by
\begin{equation*}
   u_j(x)= \left\{
\begin{array}{l}
  e^{\frac{1}{|x-x_j|^2-\eps^2}},~~~\text{if}~~x\in B_\eps(x_j) \\
  0~~~\text{if}~~x\in \Omega\backslash B_\eps(x_j).
\end{array}
\right.
\end{equation*}
These functions $u_j,$ $j=1,2$, are compactly supported and smooth (as mollifiers), so they belong to the Sobolev space $X_{\sigma,\theta}.$ Setting $u(x)=u_1(x)-u_2(x)$ for all $x\in\Omega,$ we have clearly that $u\in X_{\sigma,\theta} $ and $u\neq 0$ since $u_1$ and $u_2$ have disjoint supports. We know that 
\begin{equation*}
   |u|^{q-2}u=\sign(u)~|u|^{q-1}.
\end{equation*}
 We then have that on $\supp(u_1),$  $|u|^{q-2}u=u^{q-1}_1$ and on $\supp(u_2),$  $|u|^{q-2}u=-u^{q-1}_2.$ 
 Thus, the symmetry in the construction of $u_1$ and $u_2$, along with their disjoint supports and equal integrals, guarantees that
$$\int_\Omega |u|^{q-2}u ~dx=\int_{B_\eps(x_1)}u_1^{q-1}dx-\int_{B_\eps(x_2)}u_2^{q-1}dx=0.$$ 
Hence, we conclude that $u\in\mathcal{C}\backslash\{0\}$ as claimed.

Now, in order to state our main result, let us introduce the notations
\begin{equation}\label{Eigen-q}
\lambda_{1}(s_2,q)=\inf_{u\in \cC_{s_2,q}\setminus\{0\}}\frac{\displaystyle\iint_{\cQ}\frac{|u(x)-u(y)|^q}{|x-y|^{N+qs_2}}\ dxdy}{\displaystyle\int_{\Omega}|u|^{q}\ dx} =\inf_{u\in D_1(q)}\iint_{\cQ}\frac{|u(x)-u(y)|^q}{|x-y|^{N+qs_2}}\ dxdy
\end{equation}
where for a given $\rho>0,$
\[
D_{\rho}(q):=\big\{ u\in \cC_{s_2,q}\setminus\{0\}:\ \int_{\Omega}|u|^{q}\ dx=\rho\big\}
\]
and
\[
\cC_{s_2,q}:=\left\{u\in X_{s_2,q}:\quad\int_{\Omega}|u|^{q-2}u\ dx=0\right\}.
\]

Our main result  writes as follows.

\begin{thm}\label{Main-1} Let $s_1, s_2\in (0,1)$  and $p,q\in (1,\infty)$ with $p\neq q$. Then the set of eigenvalue  of
problem (\ref{Eq1}) is precisely $\{0\}\cup (\lambda_{1}(s_2,q),\infty)$, where $\lambda_1(s_2,q)$ is the number defined  in \eqref{Eigen-q}.
\end{thm}

Note that we assume $p\neq q$ because in the  situation $p=q$ and $s_1=s_2=s$ problem \eqref{Eq1} reduces to the eigenvalue problem considered in \cite{ML19}  and is not relevant for our discussion here. In the situation  $p=q$ and $s_1\neq s_2$, which is a purely nonlocal situation, problem \eqref{Eq1} has spectrum given by the full interval $[0,\infty)$, with the first positive eigenvalue characterized by
\[
\lambda_{1}(s_1,s_2,p)=\inf_{u\in\cC\setminus\{0\}}\frac{[u]_{s_1,p}+[u]_{s_2,p}}{\int_{\Omega}|u|^p\ dx}.
\]
We point out that with $\sigma=\max\{s_1,s_2\}$, $\theta=\max\{s_1,s_2\}$ and $p\ne q$, the Sobolev space  $X_{\sigma,\theta}$ is the possible larger space for the study  of solutions to problem \eqref{Eq1}. Unlike the local problem, where solutions is sought in the Sobolev space $W^{1,\theta}(\Omega)$, and where only two distinct complementary cases are considered with respect to parameters $p$ and $q$, with $p\neq q$, various distinguish possibilities occur in  problem \eqref{Eq1} with regard to parameters $\{s_1,s_2\}$ and $\{p,q\}$. This is due to the nonlocal feature of the operators appearing in  problem \eqref{Eq1}. We  group all the possibilities in four  distinguish cases  as follows:
\begin{itemize}
    \item[$(P_1)$] For $0<s_2< s_1< 1$ and $1<q<p<\infty$, we choose the solutions space as $X_{s_1,p}$.
    \item[$(P_2)$]For $0<s_1< s_2\le 1$ and $1<p<q<\infty$, we choose the solutions space as $X_{s_2,q}$.
    \item[$(P_3)$] For  $0<s_1\le s_2< 1$ and $1<q<p<\infty$ (crossing), we choose the solutions space as $X_{s_2,p}$ endowed with
the norm $\|u\|_{s_2,p}=[u]_{s_2,p}+\|u\|_{L^p(\Omega)}$.
    \item[$(P_4)$] For  $0<s_2\le  s_1< 1$ and $1<p<q<\infty$ (crossing), we choose the solutions space as $X_{s_1,q}$ endowed with
the norm $\|u\|_{s_1,q}=[u]_{s_1,q}+\|u\|_{L^q(\Omega)}$.
\end{itemize}
Note that for $p\neq q$, problem \eqref{pq-Laplace} is easily cover here in cases $(P_2)$ and $(P_3)$ when $s_1=s_2=1$. Moreover, in  cases $(P_2)$ and  $(P_4)$, specially when $s_2=1$ and $s_1=s\in(0,1)$, we cover the following eigenvalue problem involving mixed local-nonlocal operators,
\begin{equation}\label{Eq-Particular}
	\begin{split}
	\quad\left\{\begin{aligned}
		 (-\Delta)_p^{s}u-\Delta_qu&= \lambda |u|^{q-2}u && \text{ in \quad $\Omega$}\\
		|\nabla u|^{q-2}\partial_{\nu}u &=0     && \text{ on } \quad \partial\Omega\\
		\cN_{p,s} u&=0    && \text{ on }\quad  \mathbb{R}^N\setminus \overline{\Omega}.
	\end{aligned}\right.
	\end{split}
	\end{equation}
Eventhough equation \eqref{Eq-Particular} looks familiar to the well-known Neumann equations governed by the mixed local-nonlocal operators $-\Delta +(-\Delta)^s$, that is, when $p=q=2$, studied by many authors in recent years, see for example \cite{DV21,BDVV22,SPV22,CEF24}, we could not find  references related to the study of solutions to \eqref{Eq-Particular}. It follows from Theorem \ref{Main-1} that for  $p\neq q$,  the set of eigenvalues  of problem (\ref{Eq-Particular}) is precisely $\{0\}\cup (\lambda_{1}(q),\infty)$, where $\lambda_1(q)$ is the first positive eigenvalue of the $q$-Laplacian $-\Delta_q$ under the Neumann boundary condition.\vspace{.3cm}

We next prove the following theorem with regard to the $L^{\infty}$-bound of  any solution of  problem \eqref{Eq1} in $\R^N$, for independent interest. Note that we do not put any assumption on parameters $s_1,s_2,p$ and $q$ since other possibilities are already covered in the literature,   see \cite{SPV22} and \cite{ML19, DRV17,AFR22}, in particular, for $p=q=2$ and $0<s_1<s_2\le 1$.
\begin{thm}\label{Inf-bound}
    Assume that $s_1,s_2\in (0,1)$ and $p,q\in (1,\infty)$. Let $u\in X_{s_1,p}\cap X_{s_2,q}$ be a solution of problem \eqref{Eq1}. Then $u\in L^{\infty}(\R^N)$ and 
    \[
    \|u\|_{L^{\infty}(\R^N)}\le 2\|u\|_{L^{\infty}(\Omega)}.
    \]
\end{thm}

We use the De Giorgi iteration method in combination with the Morrey-Sobolev embedding to prove Theorem \ref{Inf-bound}. The idea of the proof of Theorem \ref{Main-1} is borrowed from the local  case.  We analyse the set of eigenvalues and the existence of eignfunctions of problem \eqref{Eq1}, depending on the possibilities listed in $(P_1)-(P_4)$ as the parameters $s_1,s_2,p$ and $q$ varied. We carefully analyse the cases $0<s_2<s_1<1<q<p<\infty$ and $0<s_1<s_2<1<p<q<\infty$. In the case $0<s_2<s_1<1<q<p<\infty$, the  functional associated to problem \eqref{Eq1} is coercive and we use the  Direct Method in the Calculus of Variations in order to find critical points of the associated energy functional. In the case $0<s_1<s_2<1<p<q<\infty$, the energy functional associated to problem \eqref{Eq1} is not  coercive. To prove that the associated energy functional has a critical point in $ X_{s_2,q} \backslash \{0\}$, we constrain the functional on the Nehari manifold and show through a series of propositions that the critical points in the Nehari manifold is in fact a solutions of problem \eqref{Eq1} in $X_{s_2,q} \backslash \{0\}$. 
The proof of the crossing cases $(P_3)$ and $(P_4)$ follows from the proof of  $(P_1)$ and $(P_2)$ respectively. This because, in  the situation $(P_3)$, the  functional associated to problem \eqref{Eq1} is coercive in $X_{s_2,p}$, and  in the situation   $(P_4)$, the energy functional associated to problem \eqref{Eq1} is not  coercive in $X_{s_1,q}$.

\vspace{.2cm}
Note that for $s_1,s_2\in (0,1)$ and $p,q\in (1,\infty)$ with $p\neq q$, our techniques  apply if one interchanges the role of the couples $(s_1,p)$ and $(s_2,q)$ in problem \eqref{Eq1}. More precisely, the results in  Theorem \ref{Main-1} and Theorem \ref{Inf-bound} hold for the problem:
  \begin{equation}\label{Eq1-p}
	\begin{split}
	\quad\left\{\begin{aligned}
		 (-\Delta)_p^{s_1}u+(-\Delta)_q^{s_2}u&= \lambda |u|^{p-2}u && \text{ in \quad $\Omega$}\\
		\cN_{p,s_1} u +\cN_{q,s_2} u&=0    && \text{ on }\quad  \mathbb{R}^N\setminus \overline{\Omega}.
	\end{aligned}\right.
	\end{split}
	\end{equation}

\vspace{.2cm}
The paper is organised as follows. In Section \ref{Preli}, we introduce the fractional Sobolev spaces in which we  work, we  review some properties and embedding results, and recall some properties related to nonlinear eigenvalue problem involving the fractional $q$-Laplacian, under nonlocal  Neumann boundary condition. Section \ref{Proof-Main} is dedicated to the proof of the main result, Theorem \ref{Main-1}. We prove nonexistence result in Section \ref{Proof-nonexistence}, existence results in Section \ref{Proof-existence} and global $L^{\infty}$-bound in section \ref{Proof-bound}.

\vspace{.2cm}
\textbf{Notations.}  In the remainder of the paper, we use the following notation. For a set $A \subset \R^N$ and $x \in \R^N$, we denote by $A^c = \R^N \setminus A$
and, if $A$ is measurable, then $|A|$ denotes its Lebesgue measure. Moreover, for given $r > 0$, let  $B_r(x) := B_r(\{x\})$ denote the ball of radius $r$ with $x$ as its center. If $x = 0$ we also write $B_r$ instead of $B_r (0)$. Moreover, if $A$ is open, we denote by $\cC_c^k(A)$ the space of function $u:\R^N \to \R$ which are $k$-times continuously differentiable and with support compactly contained in $A$. For $u:A\to \R$, we denote $u^+=\max\{u,0\}$ as the positive part of $u$ and $u^-=-\min\{u,0\}$  as the negative part of $u$, so that $u=u^+-u^-.$

\vspace{.2cm}
\textbf{Acknowledgments.}
The second author is supported by the Agence Nationale de la Recherche, Project TRECOS,
under grant ANR-20-CE40-0009.
\section{ Preliminaries}\label{Preli}

In this section, we introduce some notations, definitions and some underlying properties
of the fractional Sobolev spaces.  We also discuss the nonlinear eigenvalue problem involving the fractional $q$-Laplace operator in an open bounded set with Lipschitz boundary under nonlocal Neumann boundary condition.
Let $\Omega$ be an open bounded set of $\R^N, ~N\ge 1.$ 
We first introduce the fractional Sobolev space $X_{s,r}$.  For $0<s<1$ and $r\in(1,+\infty)$, the  space $X_{s,r}$ is defined by
\[
X_{s,r}:= \big \{ u:\R^N\to \R \quad \text{ measurable : }  \quad \|u\|_{{s,r}}<\infty\big\}
\]
where $\|u\|_{{s,r}}:= \left( \|u\|^r_{L^r(\Omega)} + [u]^r_{s,r}\right)^{1/r}~~\text{and}$
\[
[u]_{s,r}:=\Big(\iint_{\cQ}\frac{|u(x)-u(y)|^r}{|x-y|^{N+sr}}\ dxdy\Big)^{\frac{1}{r}}
\]
is the so  called  fractional  Gagliardo seminorm.   
The space $X_{s,r}$ is a reflexive  Banach space when endowed with the above norm $\|\cdot\|_{{s,r}}$ (see \cite[Proposition 2.2]{ML19}). We refer the reader to \cite{Valdinoci,L23} for further references and comprehensive examination of the fractional Sobolev Space and its properties.
We also  recall  the fractional Sobolev space $W^{s,r}(\Omega)$,  defined for $0<s<1$ and $r\in(1,+\infty)$ by
$$W^{s,r}(\Omega)=\Big\{u\in L^r(\Omega):~~ [u]^r_{s,r,\Omega}:=\iint_{\Omega\times\Omega}\frac{|u(x)-u(y)|^r}{|x-y|^{N+sr}}\ dxdy<\infty\Big\},$$ 
which,   endowed with the norm 
$\|u\|_{W^{s,r}(\Omega)}:=\left(\|u\|^r_{L^r(\Omega)}+[u]^r_{s,r,\Omega}\right)^{1/r},$ 
  is also a reflexive Banach space. 
Since the following obvious inequality holds between the usual $W^{s,r}(\Omega)$-norm and $X_{s,r}$-norm
\[
\|u\|_{W^{s,r}(\Omega)}\le C\|u\|_{{s,r}},
\]
as a consequence, $X_{s,r}\subset W^{s_2,q}(\Omega)$ and, with $\Omega$ being of finite volume, the standard fractional compact embeddings
$$\text{   $X_{r,s}\hookrightarrow L^r(\Omega)$ }$$ 
 holds for all $1<r<\infty$ and $s\in (0,1)$, see  for example \cite{Valdinoci}. Owing to \cite[Remark 2.3]{ML19}, we have that the space $X_{s,r}$ is embedded in $L^r(B(0,R))$ for every $R>0$, that is,
 $$\text{$X_{r,s}\hookrightarrow L^r(B_R(0))$~~ for every~~ $R>0$,}$$ 
 and the space $\cC^{\infty}_c(\R^N)$ is dense in $X_{s,r}$ (see \cite[Theorem 2.4]{Valdinoci}). We have  the following result from \cite[Lemma A.1]{BLP14}. 

\begin{lemma}\label{sup-semi}
    Let $1 \le  p < \infty$ and $0 < s < 1$. For every $u\in \cC_c^{\infty}(\R^N)$ it holds
    \[
    \sup_{|z|>0}\int_{\R^N}\frac{|u(x+z)-u(x)|^p}{|z|^{sp}}\ dx\le C[u]^p_{s,p}
    \]
    for some constant $C:=C(N,s,p)$.
\end{lemma}

We have  following   lemma.
\begin{lemma}\label{Embedding} 
  Let $\Omega\subset \R^N$ be a bounded set,  $0<s_2\le s_1<1< q\le p< \infty$. Then, it holds that 
  \[
  X_{s_1,p}\subset  X_{s_2,q}.\qquad 
  \]
Moreover, there exists a constant $C:=C(N, s_1, s_2, p,q, \Omega)>0$ such that 
\[
\|u\|_{s_2,q}\le C\|u\|_{s_1,p}\quad \text{for all } \quad u\in X_{s_1,p}.
\] 
As consequence, for  $0<s_2<s_1<1<q<p<\infty$  or   $0<s_1<s_2<1<p<q<\infty$, the  embedding 
\[\text{
$X_{s_1,p} \hookrightarrow X_{s_2,q}$\qquad or \qquad $X_{s_2,q} \hookrightarrow X_{s_1,p}$}
\]
is continuous.
\end{lemma}

\begin{proof} Since $\cC^{\infty}_c(\R^N)$ is dense in $X_{s,r}$,  it is enough to prove the lemma for all $u\in \cC^{\infty}_c(\R^N)$. Note  that the space $X_{s,r}$ is embedded in $L^r(B(0,R))$ for every $R>0$ (see \cite[Remark 2.3]{ML19}). Now,  Let then  $R>0$ be such that $\supp\ u\subset B(0,R)$. We have 
\begin{align*}
    [u]^q_{s_2,q}\le \int_{\R^N}\int_{\R^N}\frac{|u(x)-u(y)|^{q}}{|x-y|^{N+s_2q}}\ dxdy&=\int_{\R^N}\int_{\R^N}\frac{|u(x+z)-u(x)|^{q}}{|z|^{N+s_2q}}\ dxdz=J_1+J_2,
\end{align*}
  where the quantities $J_1$ and $J_2$ are given by
  \[
  J_1:= \int_{\{|z|<R\}}\int_{\R^N}\frac{|u(x+z)-u(x)|^{q}}{|z|^{N+s_2q}}\ dxdz,\qquad J_2:=\int_{\{|z|\ge R\}}\int_{\R^N}\frac{|u(x+z)-u(x)|^{q}}{|z|^{N+s_2q}}\ dxdz.
  \]
 We write $J_1$ as
 \begin{align*}
     J_1&=\int_{\{|z|<R\}}\frac{1}{|z|^{N-(s_1-s_2)q}}\int_{\R^N}\frac{|u(x+z)-u(x)|^{q}}{|z|^{s_1q}}\ dxdz.
 \end{align*}
 Since $u(x+z)-u(x)$ has compact support, by H\"older inequality with exponents $({p}/{q}, {p}/{(p-q)})$, we use the result in Lemma \ref{sup-semi} to get for some constant $C:=C(N,s_1,p,q,R)$,
 \begin{align*}
     \int_{\R^N}\frac{|u(x+z)-u(x)|^{q}}{|z|^{s_1q}}\ dx&= \int_{B(0,5R)}\left(\frac{|u(x+z)-u(x)|}{|z|^{s_1}}\right)^q\ dx\\
     &\le |B(0,5R)|^{\frac{p-q}{p}}\left(\int_{\R^N}\frac{|u(x+z)-u(x)|^{p}}{|z|^{s_1p}}\ dx\right)^{\frac{q}{p}}\\
     &\le |B(0,5R)|^{\frac{p-q}{p}}\sup_{|z|>0}\left(\int_{\R^N}\frac{|u(x+z)-u(x)|^{p}}{|z|^{s_1p}}\ dx\right)^{\frac{q}{p}}\\
     &\le C[u]^q_{s_1,p}.
 \end{align*}
  Hence, we estimate $J_1$ as follows, for some constant $C:=C(N,s_1,p,q,R)$, 
 \begin{align*}
     J_1&=\int_{\{|z|<R\}}\frac{1}{|z|^{N-(s_1-s_2)q}}\int_{\R^N}\frac{|u(x+z)-u(x)|^{q}}{|z|^{s_1q}}\ dxdz\\
     &\le C[u]^q_{s_1,p}\int_{\{|z|<R\}}\frac{1}{|z|^{N-(s_1-s_2)q}}\ dz\le C[u]^q_{s_1,p}.
 \end{align*}
 Next, using H\"older's inequality, we obtain for some constant $C:=C(N,s_1,p,q)$,
 \begin{align*}
     J_2&\le C(q)\int_{\{|z|\ge R\}}\int_{\R^N}\frac{|u(x)|^{q}}{|z|^{N+s_2q}}\ dxdz\\
     &\le C(q)\int_{\{|z|\ge R\}}\frac{1}{|z|^{N+s_2q}}\int_{B(0,5R)}|u(x)|^{q}\ dxdz\\
     &\le C(N,q,s_2)\int_{B(0,5R)}|u(x)|^{q}\ dx\le C\|u\|^q_{L^p(B(0,5R))}.
 \end{align*}
 Now, using the embedding $X_{s_1,p}\subset L^p(B(0,R))$, we get 
 \[
 \begin{split}
 [u]^q_{s_2,q}=\iint_{\cQ}\frac{|u(x)-u(y)|^{q}}{|x-y|^{N+s_2q}}\ dxdy&\le C(N,p,q,s_1,s_2)\left( \|u\|^q_{L^p(B(0,R))}+ [u]^q_{s_1,p}\right)\\
 &\le C [u]^q_{s_1,p}.
 \end{split}
 \]
 for some constant $C:=C(N,s_1, s_2,p,q)$. This completes the proof Lemma \ref{Embedding}.
\end{proof}
It follows from Lemma \ref{Embedding} that $X_{s_1,p}\cap X_{s_2,q}\neq\emptyset$ for $0<s_1<s_2<1<q\le p<\infty$ or $0<s_2<s_1<1<p\le q<\infty$. 
Note that  if $0<s_1<s_2<1<p\le q<\infty$, then $X_{s_1,p}\cap X_{s_2,q}\equiv X_{s_2,q}$ thanks to Lemma \ref{Embedding}. Consequently, we have that $\cC\equiv C_{s_2,q}$. Otherwise,  $\cC$ is a proper subset of $\cC_{s_2,q}$, that is when $0<s_2<s_1<1<q\le p<\infty$. 
\begin{remark}\label{Re1}
    We recall the following known fact, See \cite{MS15}. If $\Omega$ is a bounded set of $\R^N$, then 
    \[
    W^{s,p}(\Omega)\not\hookrightarrow W^{s,q}(\Omega)~~~\text{ for any }~~ 0 < s < 1 \le q < p \le  \infty.
    \]
\end{remark}\noindent
We point out that in the crossing cases (see situation $(P_3)$ and $(P_4)$), we have 
\[
  X_{s_1,p}\cap X_{s_2,q} \subsetneqq  X_{s_2,p}~~~~~\text{ and }~~~~~~ X_{s_1,p}\cap X_{s_2,q} \subsetneqq X_{s_1,q},
\]
as we shall prove in  Lemma \ref{Embedding2} below, and the space $X_{s_1,p}\cap X_{s_2,q}$ is not a suitable space for the study of solutions of problem \eqref{Eq1} in the situation $(P_3)$ and $(P_4)$. We have the following result from \cite[Page 70]{L23}.
\begin{lemma}\label{Decreas}
    Let $1\le q<p<\infty$ and $s\in (0,1)$. Let  $f:[0,\infty)\to [0,\infty)$ be an increasing function. Then, there exists a constant $C:=C(s,q)>0$ independent of $p$ such that 
    \[
    \left(\int_0^\infty (f(x))^p\frac{dx}{x^{1+sp}}\right)^{\frac{1}{p}}\le C \left(\int_0^\infty (f(x))^q\frac{dx}{x^{1+sq}}\right)^{\frac{1}{q}}.
    \]
\end{lemma}
\begin{lemma}\label{Embedding2}
Let $u\in X_{s_1,p}\cap X_{s_2,q}$ with $0<s_1<s_2<1$ and $1<q<p<\infty$. Then, there exists a constant $C:=C(s_1,s_2,p,q)>0$ such that
    \begin{equation}\label{Equiv-norm2}
[u]^p_{s_2,p}\le   C([u]^p_{s_1,p}+ [u]^p_{s_2,q}).
\end{equation}
\end{lemma}
\begin{proof} We will use the equivalent representation of the semi-norm, see \cite[Theorem 1.53]{L23}. For $0<s<1$ and $1 \le p<\infty,$ there exists two constants $C_1$ and $c_2$ such that 
\[
 C_1\sum_{j=1}^\infty\int_0^\infty \sup_{0\le t\le h}\|\Delta_{t,j}u\|^p_{L^p(\R^N)}\frac{dh}{h^{1+sp}}\le [u]^p_{s,p}\le  C_2\sum_{j=1}^\infty\int_0^\infty \sup_{0\le t\le h}\|\Delta_{t,j}u\|^p_{L^p(\R^N)}\frac{dh}{h^{1+sp}},
\]
where, denoting  $e_j$ as the $j^{th}$ vector of the canonical basis, we use the notations,
$$\Delta_{t,j}u(x)=u(x+te_j)-u(x) ~~~\text{ and }~~~~\Delta_{h}u(x)=u(x+h)-u(x).$$ 
We first note that  for $j=1,\cdots,N$, 
it follows from \cite[Lemma 7.36]{L23} that, there exists some positive  constant   $C:=C(N,p,q)$ such that
\begin{equation}\label{INQ1}
\begin{split}
   \int_0^\infty \sup_{0\le t\le h}\|\Delta_{t,j}u\|^q_{L^p(\R^N)}\frac{dh}{h^{1+s_2q}}
   &\le C  \int_0^\infty\|\Delta_{h,j}u\|^q_{L^p(\R^N)} \frac{dh}{h^{1+s_2q}}.
   \end{split}
\end{equation}
Now,  using  step $4$ in the proof of Theorem 7.35 in \cite{L23}, we get
\begin{equation}\label{INQ2}
    \left(\int_0^\infty \|\Delta_{h,j}u\|^q_{L^p(\R^N)}\frac{dh}{h^{1+s_2q}}\right)^{\frac{1}{q}}\le C\sum_{i=1}^N \left(\int_0^\infty \|\Delta_{h,i}u\|^q_{L^q(\R^N)}\frac{dh}{h^{1+s_2q}}\right)^{\frac{1}{q}}
\end{equation}
Since it holds that   $ \|\Delta_{h,j}u\|^p_{L^p(\R^N)}\le C \sup_{0\le t\le h}\|\Delta_{t,j}u\|^p_{L^p(\R^N)}$ for some constant $C>0$, and  the function $f(h)= \sup_{0\le t\le h}\|\Delta_{t,j}u\|^p_{L^p(\R^N)}$ is increasing, it follows from Lemma \ref{Decreas} and then  the inequality in \eqref{INQ1}  that 
\begin{align*}
\left(\int_0^\infty \|\Delta_{h,j}u\|^p_{L^p(\R^N)}\frac{dh}{h^{1+s_2p}}\right)^{\frac{1}{p}}&\le C\left(\int_0^\infty \sup_{0\le t\le h}\|\Delta_{t,j}u\|^p_{L^p(\R^N)}\frac{dh}{h^{1+s_2p}}\right)^{\frac{1}{p}}\\
&\le C\left(\int_0^\infty\sup_{0\le t\le h}\|\Delta_{t,j}u\|^q_{L^p(\R^N)}\frac{dh}{h^{1+s_2q}}\right)^{\frac{1}{q}}\\
&\le  C\left(\int_0^\infty \|\Delta_{h,j}u\|^q_{L^q(\R^N)}\frac{dh}{h^{1+s_2q}}\right)^{\frac{1}{q}}.
\end{align*}
We then Combine this inequality with \eqref{INQ2}, to obtain
\begin{equation}\label{INQ3}
     \left(\int_0^\infty \|\Delta_{h,j}u\|^p_{L^p(\R^N)}\frac{dh}{h^{1+s_2p}}\right)^{\frac{1}{p}}\le C\sum_{i=1}^N \left(\int_0^\infty \|\Delta_{h,i}u\|^q_{L^q(\R^N)}\frac{dh}{h^{1+s_2q}}\right)^{\frac{1}{q}}.
\end{equation}
We can conclude from \eqref{INQ3} that there exists a constant $C>0$ such that
\begin{equation}\label{INQ4}
    [u]_{s_2,p}\le C[u]_{s_2,q}.
\end{equation}
Next, in other to get \eqref{Equiv-norm2} we split the integral in the quantity $[u]^p_{s_2,p}$ as follows,
   \begin{align*}
       [u]^p_{s_2,p}= \int_{\{|z|<1\}} \frac{|u(x+z)-u(x)|^p}{|z|^{N+s_2p}}\ dxdy + \int_{\{|z|\ge 1\}} \frac{|u(x+z)-u(x)|^p}{|z|^{N+s_2p}}\ dxdy= J_1+J_2
   \end{align*} 
    where the quantities $J_1$ and $J_2$ are given by
  \[
  J_1:= \int_{\{|z|<1\}}\int_{\R^N}\frac{|u(x+z)-u(x)|^{p}}{|z|^{N+s_2p}}\ dxdz,\qquad J_2:=\int_{\{|z|\ge 1\}}\int_{\R^N}\frac{|u(x+z)-u(x)|^{p}}{|z|^{N+s_2p}}\ dxdz.
  \]
  For $J_2$, since $ps_1\le s_2p$ and $|z|\ge 1$,  $J_2$ is bounded from above by  $[u]^p_{s_1,p}$, that is,
 \[
  J_2\le \int_{\{|z|\ge 1\}}\int_{\R^N}\frac{|u(x+z)-u(x)|^{p}}{|z|^{N+s_1p}}\ dxdz\le C [u]^p_{s_1,p}.
 \]  
Now, using  spherical coordinates  and  computations  \eqref{INQ1}, \eqref{INQ2} and \eqref{INQ3},  we estimate $J_1$ as
 \begin{align*}
     J_1&=\int_{\{|z|<1\}}\int_{\R^N}\frac{|u(x+z)-u(x)|^{p}}{|z|^{N+s_2p}}\ dxdz=\int_{\{|z|<1\}}\int_{\R^N}\frac{|\Delta_hu(x)|^{p}}{|z|^{N+s_2p}}\ dxdz\\
     &=\int_{\{|z|<1\}}\frac{\|\Delta_hu\|_{L^p(\R^N)}^{p}}{|z|^{N+s_2p}}\ dz
     \le C\sum_{j=1}^N\int_{0}^1\sup_{0\le t\le h}\|\Delta_{t,j}u\|^{p}_{L^p(\R^N)} \frac{dh}{h^{1+s_2p}}\\
     &\le C\left(\sum_{j=1}^N \int_{0}^{\infty}\|\Delta_{h,j}u(x)\|^{q}_{L^q(\R^N)}\frac{dh}{h^{1+s_2q}}\right)^{\frac{p}{q}}\le C[u]^p_{s_2,q}.
 \end{align*}
 The proof of Lemma \ref{Embedding2}  is completed.
\end{proof}

We will use the following notations throughout the paper. For $u,v\in X_{s,r}$, we  introduce the functional $(u,v)\mapsto \cE_{r,s}(u,v)$, which  is the bilinear form associated to $(-\Delta)^s_r$~  defined by
$$\cE_{s,r}(u,v):= \iint_{\cQ}\frac{|u(x)-u(y)|^{r-2}(u(x)-u(y))(v(x)-v(y))}{|x-y|^{N+sr}}\ dxdy.$$
 
Next, we recall some spectral properties of the following  nonlinear eigenvalues problem involving  the fractional $q$-Laplacian $(-\Delta)^{s_2}_q$, subject to homogeneous nonlocal $q$-Neumann boundary condition, 
\begin{equation}\label{q-Neumann-eigen}
	\begin{split}
	\quad\left\{\begin{aligned}
(-\Delta)^{s_2}_q u&= \lambda  |u|^{q-2}u && \text{ in \quad $\Omega$}\\
		\cN_{q,s_2}u &=0     && \text{ on } \quad \R^N\setminus\overline\Omega.
	\end{aligned}\right.
	\end{split}
	\end{equation}
It has been shown, see \cite[Proposition 3.2]{ML19} and  \cite{DS15,DRV17} that for every $s_2\in (0,1)$ and $q\in (1,\infty)$, by means of the cohomological index and min-max theory respectively, there exists a non-decreasing sequence $\{\lambda_k(s_2,q)\}_k$ of variational eigenvalues  of $(-\Delta)_q^{s_2}$  in $\Omega$  satisfying
 $$\text{$\lambda_k(s_2,q)\to \infty$\quad as\quad $k\to \infty$}.$$
Furthermore,  the variational characterization of the $k$-eigenvalue $\lambda_{k}(s_2,q)$ is defined by the following min-max formula, 
\begin{equation}\label{min-max}
\lambda_k(s_2,q)=\displaystyle{\inf_{A\in\Sigma_k}\sup_{u\in A}\iint_{\cQ}\frac{|u(x)-u(y)|^{q}}{|x-y|^{N+s_2q}}\ dxdy}
.\end{equation}
where $\Sigma_k$ is defined by
$$
\Sigma_k=\{A\subset\Sigma\cap D(s_2,q,1):~\gamma(A)\geq k\},
$$
 the set $\Sigma,$ being the class of closed symmetric  subsets of $ X_{s_2,q}\backslash\{0\},$ i.e.,
$$
\Sigma=\{A\subset X_{s_2,q}\backslash\{0\}:~~A~~\text{closed},~~ A=-A\}
$$ 
and
\[
D( s_2, q,\rho):=\big\{u\in X_{s_2,q}\, : \, \int_{\Omega}|u|^q\ dx=\rho\big\}.
\]
For $A\in\Sigma,$ we define    
$$\gamma(A)=\inf\{k\in\mathbb{N}~:~\exists\varphi\in C(A, \mathbb{R}^k\backslash\{0\}),~\varphi(-x)=-\varphi(x)\}.$$  If such $\gamma(A)$ does not exist, we  define $\gamma(A)=+\infty.$ The number $\gamma(A)\in\mathbb{N}\cup\{+\infty\}$ is called the {\it Krasnoselskii's genus} of $A$, see \cite{AM07,FR78}.

We  next present an alternative variational characterization of the first nonzero eigenvalue $\lambda_1(s_2,q):=\lambda_{1}(s_2,q,\Omega)$ of problem \eqref{q-Neumann-eigen}. We recall that $u\in X_{s_2,q}$ is an eigenfunction of problem \eqref{q-Neumann-eigen} if 
\begin{equation}\label{Def-q}
    \cE_{s_2,q}(u,v)=\lambda \int_{\Omega}|u|^{q-2}uv\ dx~~~~~~\text{ for all} ~~v\in X_{s_2,q}.
\end{equation}

Clearly, $\lambda=0$ is an eigenvalue of \eqref{q-Neumann-eigen} with  eigenfunctions being constants. Moreover, it is isolate and simple. More precisely, we have the following proposition
\begin{prop}
    $\lambda=0$ is the first eigenvalue of problem \eqref{q-Neumann-eigen} and it is isolate and simple.
\end{prop}
\begin{proof}
    Note first that $\lambda<0$ cannot be an eigenvalue of problem \eqref{q-Neumann-eigen}. Indeed, if $\lambda$ is an eigenvalue of problem \eqref{q-Neumann-eigen}, taking $v=u$ in \eqref{Def-q} with $\|u\|_{L^q(\Omega)}=1,$ we get 
    \[
0\le [u]^q_{s_2,q}=\lambda<0,
    \]
   which is a contradiction. The simplicity of $\lambda=0$ follows from the the fact that,
   \[
   0=\inf_{u\in X_{s_2,q}\setminus\{0\}}\frac{[u]^q_{s_2,q}}{\|u\|^q_{L^q(\Omega)}}.
   \]
   Now, suppose by contradiction that $\lambda=0$ is not isolated. Then, we can find a sequence of positive eigenvalues $(\lambda_n)_n$ of  \eqref{q-Neumann-eigen}, with associated normalized  eigenfunctions $(u_n)_n$, $\|u_n\|_{L^q(\Omega)}=1,$ such that $\lambda_n\to 0$ as $n\to \infty.$ It follows that 
   \[
   \lambda_n=\frac{[u_n]^q_{s_2,q}}{\|u_n\|^q_{L^q(\Omega)}} =[u_n]^q_{s_2,q}\to 0~~~~\text{as}~~~n  \to \infty.
   \]
   Then, the sequence $(u_n)_n$ is bounded in $X_{s_2,q}$. Therefore,  passing to a
suitable subsequence that we still labeled $(u_n)_n$ , we may assume that
\[
\text{ $u_n \rightharpoonup u$ ~~~in $X_{s_2,q}$, ~~~$u_n\to u$ ~in ~$L^q(\Omega)$ ~~~and ~~~$u_n\to u$ ~a.e ~in ~$\Omega$. }
\]
We have that $\|u\|_{L^q(\Omega)}=1$ and that $[u]^q_{s_2,q}= 0$, so that $u$ is a constant, and precisely
\begin{equation}\label{Const}
u(x)= \pm\frac{1}{|\Omega|^{\frac{1}{q}}}.
\end{equation}
Moreover, taking $v=1$ in \eqref{Def-q}, we get that
\[
0=\int_{\Omega}|u_n(x)|^{q-2}u_n(x)\ dx\to 
\int_{\Omega}|u(x)|^{q-2}u(x)\ dx\neq 0.
\]
This is a contradiction with \eqref{Const} and the proof is completed.
\end{proof}
Next, suppose that $\lambda>0$ is an eigenvalue of \eqref{q-Neumann-eigen} with corresponding eigenfunction $u_{\lambda}\in X_{s_2,q}$. Then, testing \eqref{q-Neumann-eigen} with $v=1$, we obtain
\[
\int_{\Omega}|u_{\lambda}(x)|^{q-2}u_{\lambda}(x)\ dx=0.
\]
Hence,  eigenfunctions corresponding to positive eigenvalues $\lambda>0$ of problem \eqref{q-Neumann-eigen} belong to the set $\cC_{s_2,q}\setminus\{0\}$. We consider the following minimization problem
\begin{equation}\label{Lambda-N}
    \lambda_{1}(s_2,q)=\inf_{D_1(p)} \iint_{\cQ}\frac{|u(x)-u(y)|^q}{|x-y|^{N+qs_2}}\ dxdy,
\end{equation}
where $$D_1(q)=\{u\in\cC_{s_2,q}:~\|u\|^q_{L^q(\Omega)}=1\}.$$ As in the classical elliptic eigenvalue problems under Neumann boundary condition, we show that $\lambda_1(s_2,q):=\lambda_{1}(s_2,q,\Omega)$ defined in \eqref{Lambda-N},  is  an equivalent variational characterization of  the first nonzero eigenvalue of problem \eqref{q-Neumann-eigen}. We have the following proposition.
\begin{prop}Let  $s\in (0,1)$ and $p\in (1,\infty)$.
    It holds that $\lambda_1(s_2,q)>0$ and there exists $u_{\lambda_1}\in D_{1}(q)$, such that $$\lambda_{1}(s_2,q)=\iint_{\cQ}\frac{|u_{\lambda_1}(x)-u_{\lambda_1}|^q}{|x-y|^{N+qs_2}}\ dxdy.$$
    Moreover, $\lambda_1(s_2,q)$ is attained in $D_{1}(q)$ and  is the first nonzero eigenvalue of \eqref{q-Neumann-eigen}, with corresponding eigenfunction $u_{\lambda_1}$.
\end{prop}
\begin{proof}
    Let $(u_n)\subset D_1(q)$ be a minimizing sequence for $\lambda_1(s_2,q)$. Then there exists a constant $C>$ such that $\|u_n\|_{{s_2,q}}\le C$ and, evidently the sequence $(u_n)$ is bounded in $X_{s_2,q}.$ There exists a function $u\in X_{s_2,q}$ such that we may assume that
    \begin{align*}
        &u_n \rightharpoonup u ~~~\text{in }~~ X_{s_2,q},\\
        &u_n\to u ~~\text{in }~~L^q(\Omega) ~~~\text{ and} ~~~u_n\to u ~\text{ a.e ~in} ~\Omega. 
    \end{align*}
    Hence, via the Lebesgue dominated convergence theorem, we have that $\|u\|_{L^q(\Omega)}=1$ and 
    \[
0=\int_{\Omega}|u_n(x)|^{q-2}u_n(x)\ dx\to 
\int_{\Omega}|u(x)|^{q-2}u(x)\ dx.
\]
This shows that $u\in D_1(q)$. On the other hand, since $u_n\rightharpoonup u$ in $X_{s_2,q}$, exploiting the weak lower semicontinuity of the
norm functional, we have
\[
[u]^q_{s_2,q}\le \liminf_{n\to \infty} [u_n]^q_{s_2,q}= \lambda_1(s_2,q).
\]
It follows by the definition of $\lambda_1(s_2,q)$ that 
\[
[u]^q_{s_2,q}=\lambda_1(s_2,q).
\]
Since $u\in X_{s_2,q}$ and $u$ is not constant, it follows that 
$$\lambda_1(s_2,q)=[u]^q_{s_2,q}>0.$$ 
Next, by the Lagrange multiplier theorem,  we can find constants $A \in \R$ and $B \in \R$, not all of them equal to zero, such that
\begin{equation}\label{Lagrange}
\begin{split}
    \cE_{s_2,q}(u,v)+A\int_{\Omega}|u|^{q-2}uv\ dx+B\frac{p-1}{p}\int_{\Omega}|u|^{q-2}v\ dx=0,
    \end{split}
\end{equation}
for all $v\in X_{s_2,q}$. Noticing that  $u\in D_1(q)$ and letting $v=B$, we get that
\[
B^2\frac{p-1}{p}\int_{\Omega}|u|^{q-2}\ dx=0,
\]
which implies that $B=0$. Equation \eqref{Lagrange} becomes 
\begin{equation}\label{Lagrange2}
\begin{split}
    \cE_{s_2,q}(u,v)+A\int_{\Omega}|u|^{q-2}uv\ dx=0
    \end{split}
\end{equation}
for all $v\in X_{s_2,q}$. Letting $v=u$ in \eqref{Lagrange2}, we obtain
\[
\cE_{s_2,q}(u,u)+A\|u\|^q_{L^q(\Omega)}=0,
\]
which implies that $A=-\lambda_1(s_2,q)$. We then conclude that $u\in D_1(q)$ solves equation \eqref{q-Neumann-eigen}, and since clearly, there is not an eigenvalue of \eqref{q-Neumann-eigen} in the interval $(0,\lambda_1(s_2,q))$, it follows that  $\lambda_1(s_2,q)$ is the first nonzero eigenvalue of \eqref{q-Neumann-eigen}.
\end{proof}

We will employ the following optimal Poincar\'e-Wirtinger type inequality.
\begin{lemma}[Poincar\'e - Wirtinger Inequality]\label{PoincareWir} If $u\in \cC_{s_2,q}$, then
    \begin{equation}\label{Poincare}
\lambda_{1}(s_2,q)\int_\Omega|u(x)|^qdx\leq \iint_{\cQ}\frac{|u(x)-u(y)|^q}{|x-y|^{N+qs_2}}dxdy.
\end{equation}

\end{lemma}

Going back to problem \eqref{Eq1}, we consider the following minimization problem
\begin{equation}\label{Eq4}
 \lambda_1(s_1, s_2, p,q):=\displaystyle\inf_{u\in  \cC\setminus\{0\}} \frac{\displaystyle\frac{[u]^p_{s_1,p}}{p}+\frac{[u]^q_{s_2,q}}{q}}{\displaystyle\frac{1}{q}\int_\Omega |u(x)|^q\ dx } . 
\end{equation}
It holds that  
\begin{equation}\label{eigenequal}
    \lambda_1(s_1, s_2, p,q)=\lambda_{1}(s_2,q).
\end{equation}
Indeed,  for all $u\in \cC$, we clearly have that 
\begin{align*}
    \lambda_1(s_2,q)\le \lambda_1(s_1, s_2, p,q),
\end{align*}
since a positive term is added. On the other hand, consider $u=\frac1t u_{\lambda_1}$, where $u_{\lambda_1}$ is an eigenfunction of of problem \eqref{q-Neumann-eigen} associated to $\lambda_1(s_2,q)$. We get
$$\lambda_1(s_1, s_2, p,q)\leq \frac{\displaystyle\frac{[u_{\lambda_1}]^p_{s_1,p}}{pt^p}+\displaystyle\frac{[u_{\lambda_1}]^q_{s_2,q}}{qt^q}}{\displaystyle\frac{1}{qt^q}\int_{\Omega}|u_{\lambda_1}|^qdx}=\frac{\displaystyle\frac{[u_{\lambda_1}]^p_{s_1,p}}{pt^{p-q}}+\displaystyle\frac{[u_{\lambda_1}]^q_{s_2,q}}{q}}{\displaystyle\frac1q\int_{\Omega}|u_{\lambda_1}|^qdx}\to \lambda_1(s_2,q) $$
as $t\to \infty$ if $p>q$ and as $t\to 0$ if $q<q.$ We conclude that \eqref{eigenequal} holds true.

\section{Proof of the main theorem }\label{Proof-Main}

In this section, we provide a description of the set of eigenvalues of problem \eqref{Eq1}, imposing the $L^q$-normalization  $\|u_k\|^q_{L^q(\Omega)} = 1$.
We show that equation \eqref{Eq1} has  a sequence of eigenvalues $\lambda_k(s_1, s_2, p,q)$, with associated eigenfunctions $u_k$.

\subsection{Non-existence results}\label{Proof-nonexistence}
\begin{lemma}[Nonexistence] Let $\lambda_{1}(s_2,q)$ be the first eigenvalue of 
\begin{equation*}
(-\Delta)^{s_2}_q u~= ~\lambda  |u|^{q-2}u ~~~~ \text{ in \quad $\Omega$},\quad\qquad
		\cN_{s_2,q}u ~=0     ~~ \text{ on } \quad \R^N\setminus\overline\Omega.
	\end{equation*}
If it holds that $\lambda\leq \lambda_{1}(s_2,q)$, then problem (\ref{Eq1}) has no nontrivial solutions.
\end{lemma}
\begin{proof}
    Suppose by contradiction that there exists $\lambda<\lambda_{1}(s_2,q)$ which is an eigenvalue of problem (\ref{Eq1}) with $u_{\lambda}\in \cC\backslash\{0\}$ being the corresponding eigenfunction. Let $v=u_{\lambda}$ in relation (\ref{weakform}), we  have
\begin{equation*}
    \iint_{\cQ}\frac{|u_{\lambda}(x)-u_{\lambda}(y)|^p}{|x-y|^{N+s_1p}}dxdy+\iint_{\cQ}\frac{|u_{\lambda}(x)-u_{\lambda}(y)|^q}{|x-y|^{N+qs_2}}dxdy=\lambda\int_{\Omega}|u_{\lambda}(x)|^qdx.
\end{equation*}
On the other hand, the Poincar\'e inequality yields.
\begin{equation}\label{Eq5}
\lambda_{1}(s_2,q)\int_\Omega|u_{\lambda}(x)|^qdx\leq \iint_{\cQ}\frac{|u_{\lambda}(x)-u_{\lambda}(y)|^q}{|x-y|^{N+qs_2}}dxdy.
\end{equation}
Subtracting  both sides of (\ref{Eq5}) by $\lambda\int_{\Omega}|u_{\lambda}(x)|^qdx$, it follows that 
\begin{align*}
0&<\left(\lambda_{1}(s_2,q)-\lambda\right)\int_\Omega|u_{\lambda}(x)|^qdx\leq \iint_{\cQ}\frac{|u_{\lambda}(x)-u_{\lambda}(y)|^q}{|x-y|^{N+qs_2}}dxdy-\lambda \int_\Omega|u_{\lambda}(x)|^q\ dx\\
&\qquad\qquad\leq \displaystyle\iint_{\cQ}\frac{|u_{\lambda}(x)-u_{\lambda}(y)|^q}{|x-y|^{N+qs_2}}dxdy+\displaystyle\iint_{\cQ}\frac{|u_{\lambda}(x)-u_{\lambda}(y)|^p}{|x-y|^{N+s_1p}}dxdy-\lambda \displaystyle\int_\Omega|u_{\lambda}(x)|^qdx=0.
\end{align*}
Since $\lambda_{1}(s_2,q)>\lambda$, this implies that $u_{\lambda}\equiv 0$ in $\Omega$.
This contradicts the fact that $ u_{\lambda}\in \cC\backslash\{0\}$. Hence, $\lambda<\lambda_{1}(s_2,q)$ is not an eigenvalue of problem (\ref{Eq1}) with $u_{\lambda}\neq 0.$ 

The proof for the case $\lambda=\lambda_{1,q}(s_2,q)$ follows from the definition of $\lambda_{1}(s_2,q)$. Indeed, taking again $v=u_{\lambda_{1}}$ in (\ref{weakform}), we have
\begin{equation*}
    \iint_{\cQ}\frac{|u_{\lambda_1}(x)-u_{\lambda_1}(y)|^p}{|x-y|^{N+s_1p}}dxdy+\iint_{\cQ}\frac{|u_{\lambda_1}(x)-u_{\lambda_1}(y)|^q}{|x-y|^{N+qs_2}}dxdy=\lambda_{1}(s_2,q)\int_{\Omega}|u_{\lambda_1}(x)|^qdx.
\end{equation*}
It follows from the definition of $\lambda_{1}(s_2,q)$ that
\begin{equation*}
    \iint_{\cQ}\frac{|u_{\lambda_1}(x)-u_{\lambda_1}(y)|^p}{|x-y|^{N+s_1p}}dxdy\le \lambda_{1}(s_2,q)\int_{\Omega}|u_{\lambda_1}(x)|^qdx-\lambda_{1}(s_2,q)\int_{\Omega}|u_{\lambda_1}(x)|^qdx=0.
\end{equation*}
Combine this with the Poincar\'e - Wirtinger inequality in Lemma \ref{PoincareWir}, it follows  that
\begin{equation*}
0\le \int_\Omega|u_{\lambda_1}(x)|^pdx\leq C\iint_{\cQ}\frac{|u_{\lambda_1}(x)-u_{\lambda_1}(y)|^p}{|x-y|^{N+ps_1}}dxdy =0\implies u_{\lambda_1}\equiv 0\quad \text{ in }\ \Omega.
\end{equation*}
This also contradicts  the fact that $ u_{\lambda_1}\in \cC\backslash\{0\}$  and shows that any eigenvalue of problem (\ref{Eq1}) satisfies $\lambda \in (\lambda_{1}(s_2,q),\infty).$
\end{proof}

\subsection{Existence results}\label{Proof-existence}

We start the discussion about the existence of eigenvalues for problem (\ref{Eq1}). We note that these eigenvalues depend on $\rho$, from the $L^q$-normalization  $\int_{\Omega}|u(x)|^q\ dx=\rho.$ We define the  energy functional $F_{\lambda}: X_{s_1,p}\cap  X_{s_2,q}\rightarrow \mathbb{R}$ associated to relation (\ref{weakform}) by
\begin{equation}\label{functional}
    F_{\lambda}(u)=\frac{1}{p}[u]^p_{s_1,p}+\frac{1}{q}[u]^q_{s_2,q}-\frac{\lambda}{q}\int_{\Omega}|u|^qdx.
\end{equation}
A standard arguments can be used to show that $F_{\lambda}\in C^1(X_{s_1,p}\cap X_{s_2,q},\mathbb{R})$ with its derivative given by $$\langle F'_{\lambda}(u),v\rangle=\cE_{s_1,p}(u,v)+\cE_{s_2,q}(u,v)-\lambda\int_{\Omega}|u|^{q-2}u~v~dx,$$ for all $v\in X_{s_1,p}\cap X_{s_2,q}$. Thus, we note that $\lambda$ is an eigenvalue of problem (\ref{Eq1}) if and only if $F_{\lambda}$ possesses a nontrivial critical point. 

We further split the discussion into two cases, whether, $0<s_2< s_1<1<q<p<\infty$ ~ or ~  $0<s_1< s_2<1<p<q<\infty$.

\subsubsection{The case: $0<s_2< s_1<1<q<p<\infty$}

In this case $X_{s_1,p}\cap X_{s_2,q}\equiv X_{s_1,p}$ and $\cC\subset \cC_{s_2,q}$. We will show that  for each $\lambda>0,$ the functional $F_{\lambda}$ defined in (\ref{functional}) is coercive.

\begin{lemma}\label{coercivite}
 For each $\lambda>0,$ the functional $F_{\lambda}$ defined in (\ref{functional}) is coercive, that is, 
 \[
 \lim_{\substack{\|u\|_{{s_1,p}}\to +\infty\\u\in\cC}}F_{\lambda}(u) =+\infty.
 \]
\end{lemma}
\begin{proof}
Note that for $u\in \cC\subset \cC_{s_2,q}$, we get $[u]^q_{s_2,q}\ge \lambda_1(s_2,q)\|u\|^q_{L^q(\Omega)}$ by the Poincar\'e - Wirtinger inequality (Lemma \ref{PoincareWir}). However, with $q<p$, it holds using H\"older inequality that $\|u\|_{L^q(\Omega)}\le |\Omega|^{\frac{p-q}{pq}}\|u\|_{L^p(\Omega)}$.  Now, using the fact that $X_{s_1,p}\subset X_{s_2,q}$ (see Lemma \ref{Embedding}), which implies that $\|u\|_{s_2,q}\le C\|u\|_{s_1,p}$ for some $C>0$, we have
\begin{align*}
   F_{\lambda}(u)&\ge  \frac{[u]^p_{s_1,p}}{p}+\frac{\lambda_1(s_2,q)}{q}\int_{\Omega}|u|^qdx-\frac{\lambda}{\lambda_1(s_2,q)}[u]^q_{s_2,q}\\
   &\ge  \frac{1}{p}\min\{1,\lambda_1(s_2,q)|\Omega|^{\frac{p-q}{q}}\}\|u\|^p_{{s_1,p}}-\frac{C}{q}\|u\|_{{s_1,p}}^q= C_1\|u\|^p_{{s_1,p}}-C_2\|u\|^q_{{s_1,p}}.
\end{align*}
Since $p>q,$ it follows that
 \[
 \lim_{\substack{\|u\|_{{s_1,p}}\to +\infty\\u\in\cC}}F_{\lambda}(u) =+\infty.
 \]
 This shows that $F_{\lambda}$ is coercive.
\end{proof}
\begin{remark}
Note that the functional $F_{\lambda}$ is not bounded below if $p<q$ and $\lambda>\lambda_1(s_2,q)$. Indeed,  for  $u=u_1,$ the first eigenfunction of problem (\ref{q-Neumann-eigen}) with $\|u_1\|_{L^q(\Omega)}^q=1,$ we have 
\begin{align*}
F_{\lambda}(tu_1)=\frac{t^p}{p}[u_1]^p_{s_1,p}-\frac{t^q}{q}(\lambda-\lambda_1(s_2,q))\rightarrow -\infty\qquad \text{as \quad $t\rightarrow+\infty.$}
\end{align*}
This case will be treated in Section \ref{case2}, on a subset of $X_{s_1,p}\cap  X_{s_2,q}$, the so-called Nehari manifold, since we cannot apply the Direct Method in the Calculus of Variations in order to find critical points for the functional $F_{\lambda}$. 
\end{remark}

\begin{thm}\label{Existence1}
Every $\lambda\in  (\lambda_1(s_2,q),+\infty)$ is an eigenvalue of problem (\ref{Eq1}).
\end{thm}
\begin{proof}
As already  mentioned above, standard arguments show that $F_{\lambda}\in C^1(X_{s_1,p},\mathbb{R})$ with its derivative given by 
$$\langle F'_{\lambda}(u),v\rangle=\cE_{s_1,p}(u,v)+\cE_{s_2,q}(u,v)-\lambda\int_{\Omega}|u|^{q-2}u~v~dx,$$
for all $v\in X_{s_1,p}\subset X_{s_2,q}.$ On the other hand, since $\cC$ is a weakly closed subset of $X_{s_1,p}$ and 
$F_{\lambda}$ is weakly lower semi-continuous on $\cC$, since $F_{\lambda}$ is a continuous convex functional. This fact and Lemma \ref{coercivite} allow one to apply the Direct Method of Calculus of Variations   to obtain the existence of global minimum point of $F_{\lambda}$. We denote by $u_0$ such a global minimum point, i.e, 
$$F_{\lambda}(u_0)=\min\limits_{u\in \cC }F_{\lambda}(u).$$ 
We observe that for $u_0=tw_1$, where $w_1$ stands for the $L^q$-normalized associated eigenfunction of $\lambda_1(s_2,q)=\lambda_1(s_1,s_2,p,q)$, we have, 
$$F_{\lambda}(u_0)=F_{\lambda}(tw_1)=\frac{t^p}{p}[w_1]^p_{s_1,p}+\frac{t^q}{q}(\lambda_1(s_2,q)-\lambda)<0$$
for $t$ small enough. So there exists $u_{\lambda}\in \cC$ such that $F_{\lambda}(u_{\lambda})<0.$ But $F_{\lambda}(u_0)\leq F_{\lambda}(u_{\lambda})<0,$ which implies that $u_0\in \cC\backslash\{0\}.$ 

We next show that 
$$\text{$\langle F'_{\lambda}(u_0),v\rangle=0$,\quad for all\quad \ $v \in X_{s_1,p},$}$$ 
under the constraint 
\[
\int_{\Omega}|u_0|^{p-2}u_0\ dx=0.
\]
By the Lagrange multiplier theorem, there exists a constant  $B\in\R$, such that
\begin{equation}\label{Lagrange3}
\langle F'_{\lambda}(u_0),v\rangle=B(p-1)\int_{\Omega}|u_0|^{p-2}v\ dx \quad \text{for all}\quad \ v \in X_{s_1,p}.
\end{equation}
Letting $ v = B$ in the above equation \eqref{Lagrange3}, we get that
\[
B^2(p-1)\int_{\Omega}|u_0|^{p-2}\ dx=0,
\]
which implies that $B = 0$. Equation \eqref{Lagrange3} becomes
\begin{equation}\label{Lag}
\langle F'_{\lambda}(u_0),v\rangle=0 \quad \text{for all}\quad \ v \in X_{s_1,p}.
\end{equation}
We  conclude that $u_0$ solves equation \eqref{Eq1} and any $\lambda\in (\lambda_1(s_2,q),+\infty)$ is an eigenvalue of \eqref{Eq1}. This concludes the proof of Theorem \ref{Existence1}.
\end{proof}

\subsubsection{The case: $0<s_1< s_2<1<p<q<\infty$ }\label{case2}

In this case, we do not have  coercivity on $X_{s_2,q}\setminus\{0\}$ for the functional $F_{\lambda}$ eventhough it belongs to $\cC^1(X_{s_2,q}\setminus\{0\};\R)$. Moreover, with the assumptions on $(s_1,p)$ and $(s_2,q)$, we have that $X_{s_1,p}\cap X_{s_2,q}\equiv X_{s_2,q}$ and $\cC\equiv\cC_{s_2,q}$.   To prove that $F_{\lambda
}$ has a critical point in $\cC\setminus\{0\}$, we constrain $F_{\lambda}$ on the Nehari set
\begin{equation*}
\begin{split}
    \mathcal{N}_{\lambda}&=\{u\in \cC_{s_2,q}\setminus\{0\},~\langle F'_{\lambda}(u),u\rangle=0\}\\
    &=\{u\in \cC_{s_2,q}\setminus\{ 0\},~[u]^p_{s_1,p}+[u]^q_{s_2,q}=\lambda\int_{\Omega}|u|^q\ dx\}.
\end{split}
\end{equation*}
Note that on $\mathcal{N}_{\lambda},$ the functional $F_{\lambda}$ reads as $$F_{\lambda}(u)=(\frac{1}{p}-\frac{1}{q})[u]^p_{s_1,p}>0.$$
This shows at once that $F_{\lambda}$ is coercive in the sense that if $u\in\mathcal{N}_{\lambda}$ satisfies $[u]_{s_1,p}\rightarrow+\infty,$ then $F_{\lambda}(u)\rightarrow +\infty.$

Next, we define the quantity
$$\mathfrak{M}=\inf\limits_{u\in\mathcal{N}_{\lambda}}F_{\lambda}(u)$$
and  show through a series of propositions that~ $\mathfrak{M} $ is attained by a function $u\in\mathcal{N}_{\lambda}$, which is a critical point of $F_{\lambda}$ considered on the whole space $\cC$ and in fact is a solution of (\ref{Eq1}).

\begin{prop}\label{nonvide}
The set $\mathcal{N}_{\lambda}$ is not empty for any $\lambda>\lambda_1(s_2,q).$
\end{prop}
\begin{proof}
Since $\lambda>\lambda_1(s_2,q)$ there exists $u\in \cC_{s_2,q}$ not identically zero such that 
$$[u]^q_{s_2,q}<\lambda\int_{\Omega}|u|^qdx.$$ 
Since $tu\in \cC_{s_2,q}$ for all $t>0$, there exists $t>0$, $tu\in \mathcal{N}_{\lambda}$.  Indeed, $tu\in \mathcal{N}_{\lambda}$ is equivalent to
$$t^p[u]^p_{s_1,p}+t^q[u]^q_{s_2,q}=t^q\lambda\int_{\Omega}|u|^q\ dx,$$ which is solved for $$t=\left(\frac{[u]^p_{s_1,p}}{\lambda\int_{\Omega}|u|^q\ dx-[u]^q_{s_2,q}}\right)^{\frac{1}{q-p}}>0.$$
This completes the proof of Proposition \ref{nonvide}.
\end{proof}

\begin{prop}\label{minborne}
Every minimizing sequence for $F_{\lambda}$ on $\mathcal{N}_{\lambda}$ is bounded in $X_{s_2,q}.$
\end{prop}
\begin{proof}
Let $\{u_n\}_{n\geq0}\subset \mathcal{N}_{\lambda}$ be a minimizing sequence of $F_{\lambda}|_{\mathcal{N}_{\lambda}}$, i.e. $F_{\lambda}(u_n)\rightarrow \mathfrak{M}=\displaystyle{\inf_{v\in\mathcal{N}_{\lambda}}}F_{\lambda}(v).$
Then
\begin{equation}\label{m1}
\lambda\int_{\Omega}|u_n|^q~dx-
[u_n]^q_{s_2,q}=[u_n]^p_{s_1,p}\rightarrow\left(\frac{1}{p}-\frac{1}{q}\right)^{-1}\mathfrak{M},~\text{as $n\rightarrow +\infty$}.
\end{equation}
Suppose on the contrary that $\{u_n\}_{n\geq0}$ is not bounded in $X_{s_2,q}$ i.e. $[u_n]^q_{s_2,q}\rightarrow +\infty$ as $n\rightarrow +\infty$. Then we have $\displaystyle{\int_{\Omega}}|u_n|^q\ dx
\rightarrow\infty$ as $n\rightarrow +\infty$, using relation (\ref{m1}). Next, we set $w_n=\frac{u_n}{\|u_n\|_{L^q(\Omega)}}.$ Since $[u_n]^q_{s_2,q}<\lambda\displaystyle{\int_{\Omega}}|u_n|^q~dx
$, we deduce that $[w_n]^q_{s_2,q}<\lambda,$ for each $n$. Hence $\{w_n\}_n$ is uniformly bounded in $X_{s_2,q}.$ Therefore, there exists $w_0\in X_{s_2,q}$ such that passing to a subsequence, that we still label $\{w_n\}_n$, we may assume that
\begin{align*}
&w_n\rightharpoonup w_0\quad \text{ in }\quad X_{s_2,q}\subset X_{s_1,p},\\
&w_n\rightarrow w_0 \quad\text{ in } \quad L^q(\Omega) \quad \text{ and }\quad w_n\rightarrow w_0 \quad\text{a.e~~ in} \quad \Omega.
\end{align*}
It follows by the Lebesgue dominated convergence theorem that  
    \[
0=\int_{\Omega}|u_n(x)|^{q-2}u_n(x)\ dx\to 
\int_{\Omega}|u(x)|^{q-2}u(x)\ dx,
\]
and hence $w_0\in \cC_{s_2,q}$.
 Since $p<q$, by  the H\"older inequality, we have that
$
\|w_n\|_{L^p(\Omega)}\le |\Omega|^{\frac{q-p}{pq}}\|w_n\|_{L^q(\Omega)}\le C 
$
and $w_n\to w_0$ in $L^p(\Omega)$ as $n\to \infty$.
Furthermore, as we know that $\lambda\|u_n\|^q_{L^q(\Omega)}-
[u_n]^q_{s_2,q}$ is bounded as $n\to \infty$ and $p<q$, it follows from (\ref{m1})  that
\[
\begin{split}
& [w_n]^p_{s_1,p}=\ \frac{\displaystyle\lambda\int_{\Omega}|u_n|^q \ dx-
[u_n]^q_{s_2,q}}{\|u_n\|^p_{L^q(\Omega)}}\to 0\quad \text{ as }\quad n\to \infty.
\end{split}
\]
Therefore, $w_0$ is a constant and belongs to $\cC$.   Consequently $w_0=0$. This contradicts the normalization of $\|w_n\|_{L^q(\Omega)}=1$. Hence the sequence $\{u_n\}_n$ is bounded in $X_{s_2,q}$. This completes the proof of Proposition \ref{minborne}.
\end{proof}

\begin{prop}\label{minpositive}
We have ~ $\mathfrak{M}=\inf\limits_{u\in\mathcal{N}_{\lambda}}F_{\lambda}(u)>0.$
\end{prop}
\begin{proof}
Assume by contradiction that $\mathfrak{M}=0$. Then, with $\{u_n\}_{n}$ as in Proposition \ref{minborne}, we have
\begin{equation}\label{e6}
0<\lambda\int_{\Omega}|u_n|^q \ dx-
[u_n]^q_{s_2,q}=[u_n]^p_{s_1,p}~\rightarrow ~0, \quad\text{as ~~~$n\rightarrow +\infty$}.
\end{equation}
By Proposition \ref{minborne}, we know that $\{u_n\}_{n}$ is bounded in $X_{s_2,q}.$ Therefore there exists $u_0\in X_{s_2,q}$ such that up to a subsequence,  $u_n \rightharpoonup u_0$ in $X_{s_2,q}\subset X_{s_1,p}$ and $u_n\rightarrow u_0$ in $L^q(\Omega).$
Thus, by  the weakly lower semicontinuity of the norm,
$$[u_0]^p_{s_1,p}\leq C\liminf\limits_{n\rightarrow +\infty}[u_n]^p_{s_1,p}=0.$$ 
 Consequently  $u_0$ is a constant and  $u_0\in\cC,$. It follows that $u_0=0$. Passing to a subsequence again, we have that $u_n \rightharpoonup 0$ in $X_{s_2,q}\subset X_{s_1,p}$ and $u_n\rightarrow 0$ in $L^q(\Omega) $ and also $u_n\rightarrow 0$ in $L^p(\Omega) $ since $p<q.$  Next,  set again $w_n=\frac{u_n}{\|u_n\|_{L^q(\Omega)}}$, we have that $\|w_n\|_{L^q(\Omega)}=1$ for all $n$ and as in Proposition \ref{minborne}, $\{w_n\}_n$ is uniformly bounded in $X_{s_2,q}.$ Hence $w_n \rightharpoonup w_0$ in $X_{s_2,q}$ and $w_n\rightarrow w_0$ in $L^q(\Omega)$ and $w_n\rightarrow w_0$ in $L^p(\Omega)$ since $p<q$.  By Poincar\'e inequality, we have
\[
\begin{split}
      [w_n]^p_{s_1,p}
 &=\frac{C }{\|u_n\|^p_{L^q(\Omega)}}\left(\lambda\int_{\Omega}|u_n|^q~dx-
[u_n]^q_{s_2,q}  \right) \\
&=C \|u_n\|^{q-p}_{L^q(\Omega)}\left(\lambda-
[w_n]^q_{s_2,q}  \right) \to 0\quad \text{ as }\quad n\to \infty.
\end{split}
\]
We deduce  taking the limit that $w_0$ is a constant and belongs to $\cC$. Therefore $w_0=0,$ which contradicts the normalization $\|w_n\|_{L^q(\Omega)}=1$. This completes the proof of Proposition \ref{minpositive}.
\end{proof}

\begin{prop}\label{minatteind}
 There exists $u\in \mathcal{N}_{\lambda}$ such that
$F_{\lambda}(u)=\mathfrak{M}.$
\end{prop}
\begin{proof}
Let $\{u_n\}_{n}\subset\mathcal{N}_{\lambda}$ be a minimizing sequence, i.e., $F_{\lambda}(u_n)\rightarrow \mathfrak{M}$ as $n\rightarrow\infty.$ Thanks to Proposition \ref{minborne}, we have that $\{u_n\}$ is bounded in $X_{s_2,q}.$ It follows that there exists $u_0\in X_{s_2,q}$ such that $u_n\rightharpoonup u_0$ in $X_{s_2,q}\subset X_{s_1,p}$ and strongly in $L^q(\Omega).$ The results in the two propositions above guarantee that $F_{\lambda}(u_0)\leq \displaystyle{\lim_{n\rightarrow\infty}}\inf F_{\lambda}(u_n)=\mathfrak{M}.$ Since for each $n$ we have $u_n\in\mathcal{N}_{\lambda}$, then
\begin{equation}\label{att1}
[u_n]^p_{s_1,p}+[u_n]^q_{s_2,q}=\lambda \int_{\Omega}|u_n|^q~dx~~~\text{for all $n$.}
\end{equation}
Assuming $u_0\equiv 0$ on $\Omega$ implies that $ \displaystyle{\int_{\Omega}}|u_n|^q~dx\rightarrow 0$ as $n\rightarrow +\infty$, and by relation $(\ref{att1})$ we obtain that $[u_n]^q_{s_2,q}\rightarrow 0$ as $n\rightarrow +\infty.$ Combining this with the fact that $u_n \rightharpoonup 0$ in $X_{s_2,q}$, we deduce that $u_n$ converges strongly to $0$ in $X_{s_2,q}$ and consequently in $X_{s_1,p}$.
Hence, we infer that \begin{eqnarray*}
\lambda\int_{\Omega}|u_n|^q~dx-
[u_n]^q_{s_2,q}=[u_n]^p_{s_1,p}\rightarrow 0, \text{as $n\rightarrow +\infty$}.
\end{eqnarray*}
Next, using similar argument as the one used in the proof of Proposition \ref{minpositive}, we will reach to a contradiction, which shows that $u_0\not\equiv 0.$ Letting $n\rightarrow \infty$ in relation (\ref{att1}), we deduce that 
\begin{eqnarray*}
[u_0]^p_{s_1,p}+[u_0]^q_{s_2,q}\leq\lambda \int_{\Omega}|u_0|^q~dx.
\end{eqnarray*}
If there is equality in the above relation then $u_0\in\mathcal{N}_{\lambda}$ and $\mathfrak{M}\leq F_{\lambda}(u_0)$. Assume by contradiction that
\begin{equation}\label{att2}
[u]^p_{s_1,p}+[u]^q_{s_2,q}<\lambda \int_{\Omega}|u|^q~dx.
\end{equation} 
Let $t>0$ be such that $tu_0\in\mathcal{N}_{\lambda},$ i.e.,$$t=\Bigg(\frac{\lambda\displaystyle{\int_{\Omega}}|u_0|^q~dx-
[u_0]^q_{s_2,q}}{[u_0]^p_{s_1,p}}\ \Bigg)^{\frac{1}{p-q}}
.$$
We note that $t\in(0,1)$ since $1<t^{p-q}$ (using (\ref{att2})). Finally, since $tu_0\in \mathcal{N}_{\lambda}$ with $t\in(0,1)$ we have 
\begin{eqnarray*}
0<\mathfrak{M}\leq F_{\lambda}(tu_0)&=&\left(\frac{1}{p}-\frac{1}{q}\right)[tu_0]^p_{s_1,p}=t^p\left(\frac{1}{p}-\frac{1}{q}\right)[u_0]^p_{s_1,p}\\
&=&t^p F_{\lambda}(u_0)\\
&\leq & t^p\lim_{n\rightarrow+\infty}\inf F_{\lambda}(u_n)=t^p \mathfrak{M}<\mathfrak{M}~~~~\text{for $t\in(0,1)$.}
\end{eqnarray*}
This is a contradiction, which assures that relation (\ref{att2}) cannot hold and consequently we have $u_0\in \mathcal{N}_{\lambda}$. Hence, $\mathfrak{M}\leq F_{\lambda}(u_0)$ and $ \mathfrak{M}= F_{\lambda}(u_0)$.
\end{proof}

\begin{thm}\label{constantsign2}
Every $\lambda\in (\lambda_1(s_2,q),+\infty)$ is an eigenvalue  of problem (\ref{Eq1}). 
\end{thm}
\begin{proof}
Let $u\in \mathcal{N}_{\lambda}$ be such that $F_{\lambda}(u)=\mathfrak{M}$ thanks to Proposition \ref{minatteind}. We show that 
$$\text{$\langle F'_{\lambda}(u), v\rangle=0$\quad  for all \quad $v\in X_{s_2,q}.$}$$ 
We recall that for $u\in \mathcal{N}_{\lambda}$, we have \begin{eqnarray*}
[u]^p_{s_1,p}+[u]^q_{s_2,q}=\lambda \int_{\Omega}|u|^q~dx.
\end{eqnarray*} 
 Since also  $u\neq 0$, we have that $[u]^q_{s_2,q}<\lambda \int_{\Omega}|u|^q~dx$. Let $v\in \cC$ be arbitrary and fixed. Then, there exists $\epsilon>0$ such that  for every $\delta$ in the interval $(-\varepsilon,\varepsilon)$, certainly the function $u+\delta v$ does not vanish identically and $[u+\delta v]^q_{s_2,q}<\lambda \int_{\Omega}|u+\delta v|^q~dx$.~ Let $t(\delta)>0$ be a function such that $t(\delta)(u+\delta v)\in  \mathcal{N}_{\lambda},$ namely 
$$t(\delta)=\Bigg(\frac{\lambda\displaystyle{\int_{\Omega}}|u+\delta v|^q~dx-
[u+\delta v]^q_{s_2,q}}{\displaystyle{[u+\delta v]^p_{s_1,p}}}\ \Bigg)^{\frac{1}{p-q}}
.$$ The function $t(\delta)$ is a composition of differentiable functions, so it is differentiable. The precise expression of $t'(\delta)$ does not matter here. Observe that $t(0)=1.$ The map $\delta\mapsto t(\delta)(u+\delta v) $
defines a curve on $\mathcal{N}_{\lambda}$ along which we evaluate $F_{\lambda}.$ Hence we define ~$\ell: (-\varepsilon,\varepsilon)\rightarrow \mathbb{R}$ as
$$\text{ $\ell(\delta)=F_{\lambda}(t(\delta)(u+\delta v)).$}$$ By construction, $\delta=0$ is a minimum point for $\ell.$ Consequently $$0=\ell'(0)=\langle F'_{\lambda}(t(0)u), t'(0)u+t(0)v\rangle=t'(0)\langle F'_{\lambda}(u),u\rangle+\langle F'_{\lambda}(u),v\rangle=\langle F'_{\lambda}(u),v\rangle.$$ Using the fact that $\langle F'_{\lambda}(u), u\rangle=0$ for $u\in \mathcal{N}_{\lambda}$, we obtained that 
$$\text{$\langle F'_{\lambda}(u),v\rangle=0$\quad  for all \quad $v\in\cC$}.$$
The proof of Theorem \ref{constantsign2} is completed.
\end{proof}

\subsubsection{The crossing cases}
In this section, we deal with the situations $(P_3)$ and ($P_4)$, that is, when    $0<s_1\le s_2\le 1$ and $1<q<p<\infty$,  and when  $0<s_2\le  s_1\le 1$ and $1<p<q<\infty$ respectively. We recall that for any fixed  $\lambda>\lambda_1(s_2,q)$ the functional associated to problem \eqref{Eq1} is 
\begin{equation}\label{b}
    F_{\lambda}(u)=\frac{1}{p}[u]^p_{s_1,p}+\frac{1}{q}[u]^q_{s_2,q}-\frac{\lambda}{q}\int_{\Omega}|u|^qdx.
\end{equation}
Note that in the situation $(P_3)$, that is, with  $0<s_1\le s_2\le 1$ and $1<q<p<\infty$, it follows from Lemma \ref{Embedding}, see also \cite[Proposition 2.1]{Valdinoci} and \cite[Corolla 7.38]{L23}, since $s_2-N/p\ge s_1-N/p$, that there exists a constant $C:=C(N,p,q,s_1,s_2)>0$ such that
\begin{equation}\label{CRembed}
\|u\|_{s_1,p}\le C\|u\|_{s_2,p}
\end{equation}
We recall that the space $X_{s_2,p}$ is endowed with the norm $\|u\|_{s_2,p}=  [u]_{s_2,p}+\|u\|_{L^p(\Omega)}$. Since $p>q$, it follows from \cite[Remark 15, P. 286]{B11} that the norm $\|\cdot\|_{s_2,p}$ is equivalent to 
\begin{equation}\label{Equiv-norm}
|||u|||_{s_2,p}=  [u]_{s_2,p}+\|u\|_{L^q(\Omega)}.
\end{equation}
Next, in the situation $(P_3)$,  we can observe that for $u=\alpha u_1$, where $u_1$ is an  eigenfunction associated to the first positive eigenvalue $\lambda_{1}(s_2,q)$ of the fractional $p$-Laplacian with $\|u_1\|_{L^q(\Omega)}=1$, substituting $u_1$ into the functional $F_\lambda$ it follows that
\begin{align*}
    F_\lambda(\alpha u_1)&= \frac{\alpha^p}{p}[u_1]^p_{s_1,p}+\frac{\lambda_{1}(s_2,q)\alpha^q}{q}-\frac{\lambda~\alpha^q}{q}=\alpha^p\Big(\frac{1}{p}[u_1]^p_{s_2,p}+\frac{\alpha^{q-p}}{q}(\lambda_1(s_2,q)-\lambda)\Big).
\end{align*}
 Therefore, $F_\lambda(\alpha u)\to \infty ~ \text{ as } ~\alpha\to\infty$ for every $u\in  X_{s_2,p}$ thanks to \eqref{CRembed}. This indicates  that the functional $F_\lambda$ is coercive for $p>q$, as we shall prove in Lemma \ref{coercive2} below.\vspace{.2cm}

In the situation $(P_4)$, that is, when  $0<s_2\le s_1\le 1$ and $1<p<q<\infty$, it follows from Lemma \ref{Embedding}, see also \cite[Proposition 2.1]{Valdinoci}, that for some constant $C:=C(N,p,q,s_1,s_2)>0$,
$$\|u\|_{s_2,q}\le C\|u\|_{s_1,q}.$$
  As above, we have with $\lambda>\lambda_{1}(s_2,q)$ with $u=\alpha u_1$,  that
\begin{align*}
    F_\lambda(\alpha u_1)&=\alpha^q\Big(\frac{\alpha^{p-q}}{p}[u_1]^p_{s_1,p}+\frac{1}{q}(\lambda_1(s_2,q)-\lambda)\Big)~\to~ -\infty ~ ~\text{ as } ~~\alpha\to\infty.
\end{align*}
This shows that the functional $F_\lambda$ is not bounded bellow if $q<p$. Hence the situation $(P_4)$ will be treated on the Nehari manifold which is a subset of  $X_{s_1,q}$. 
\begin{lemma} \label{coercive2}
Under the assumption $(P_3)$, the functional $J_\lambda$ is coercive, that is 
\[
\lim_{\substack{\|u\|_{s_2,p}\to +\infty\\u\in \cC}}\left(\frac{1}{p}[u]^p_{s_1,p}+\frac{1}{q}[u]^q_{s_2,q}-\frac{\lambda}{q}\int_{\Omega}|u|^qdx\right) =+\infty.
\]
\end{lemma}
\begin{proof} We shall use the equivalence of the  norms $\|u\|_{s_2,p}$ and $\||u\||_{s_2,p}$. Hence $\|u\|_{s_2,p}\to +\infty$ if and only if  $|||u|||_{s_2,p}\to +\infty$. Since $p>q$, we have that
\[
\frac{1}{p}[u]^p_{s_1,p}+\frac{1}{q}[u]^q_{s_2,q}\ge \frac{1}{p}\big([u]^p_{s_1,p}+[u]^q_{s_2,q}\big).
\]
It follows from Lemma \ref{Embedding2} that, 
\[
\lim_{\substack{\|u\|_{s_2,p}\to +\infty\\u\in\cC}}\left(\frac{1}{p}[u]^p_{s_1,p}+\frac{1}{q}[u]^q_{s_2,q}\right)\to +\infty.
\]
Next,   we use Lemma \ref{Embedding2} and Lemma \ref{PoincareWir}, combined with H\"older inequality and \eqref{b}, to obtain the following,
\begin{align*}
     F_{\lambda}(u)&=\frac{1}{p}[u]^p_{s_1,p}+\frac{1}{q}[u]^q_{s_2,q}-\frac{\lambda}{q}\int_{\Omega}|u|^qdx\\
    &\geq \frac1p\big([u]^p_{s_1,p}+[u]^q_{s_2,q}\big)-\frac{\lambda}{q}\int_{\Omega}|u|^qdx\\
    &\geq \frac1p\big([u]^p_{s_1,p}+[u]^q_{s_2,q}\big)-\frac{\lambda |\Omega|^{(p-q)/p}}{ q\lambda_1(s_2,q)}|||u|||^q_{s_2,p}.
\end{align*}
Since $p>q>1,$ we infer that the term in the right-hand side of the above inequality blows up as 
$\|u\|_{s_2,p}\to +\infty.$
\end{proof}
\vspace{.2cm}
We now complete the proof of Theorem \ref{Main-1}
\begin{proof}[Proof of Theorem \ref{Main-1}] The proof of Theorem \ref{Main-1} in the situation $(P_1)$ follows from Theorem \ref{Existence1} and the proof of  Theorem \ref{Main-1} in the situation $(P_2)$ follows from Theorem \ref{constantsign2}. 

In the situation $(P_3)$,  $F_\lambda$ is coercive thanks to Lemma \ref{coercive2}. We use Direct Method in the Calculus of Variations in order to find critical points of   $F_\lambda$ as the situation $(P_1)$. 

In the situation $(P_4)$, since the functional $F_\lambda$ is not coercive,  we constrain the functional $F_\lambda$ on the Nehari manifold $\cN_\lambda$  and follow the same steps of the proof as in  the situation $(P_2).$
\end{proof}

\section{ Global  $L^{\infty}$-bound} \label{Proof-bound}

We use the De Giorgi iteration method in combination with the Morrey-Sobolev embedding to prove Theorem \ref{Inf-bound}. Moreover, we  will need the following lemma.

\begin{lemma}\label{Giusti}\cite[Lemma 7.1, P.220]{Giusti03}.
Let $(U_n)_{n}$ be a sequence of non-negative  real number such that  
\[
U_{n+1}\le Cb^{n}U_n^{1+\alpha}\qquad \text{for all }\quad n,
\]
where $C>0$, $b>1$ and $\alpha>0$. Assume that $U_0\le C^{-\frac{1}{\alpha}}b^{-\frac{1}{\alpha^2}}$. Then~ $U_n \to 0$ as $n\to \infty$.
\end{lemma}

 We now give the proof of Theorem \ref{Inf-bound}.
 
\begin{proof}[Proof of Theorem \ref{Inf-bound}] We first prove that $u$ is bounded in $\Omega$. Note that if $qs_2>N$ the conclusion is a consequence of Morrey embedding. We next assume that $N>qs_2$ and set $q_{s_2}^*=2N/(N-qs_2)$, the critical Sobolev exponent. We will only  prove that $u^+$ is bounded in $\Omega$, and since both $u^{\pm}$ are solutions,  the bound for the negative part $u^-$ will follow in the same way by replacing $u$ with $-u$ in our computation. It is enough to prove that
\begin{equation}
    \|u\|_{L^{\infty}(\Omega)}\le 1 ~~~~~~~~\text{ if }~~~~~\|u\|_{L^q(\Omega)}\le\delta,
\end{equation}
where $\delta>0$ is a small parameter to be determined later. Next, for any $n\in \N$, we set $C_n: = 1-2^{-n}$ and consider the family of truncated functions
\[
w_{n}:= \left(u-C_n\right)^+, ~~~~~~~~U_n:= \int_{\Omega}|w_{n}(x)|^q\ dx.
\]
Since $C_{n+1}>C_n$, for all $n$, the sequence $\{C_n\}_{n}$ is increasing  and  consequently 
$$\text{$w_{n+1}\le w_{n}$ ~~a.e.~~ in $\R^N$.}$$ 
Moreover, since $\Omega$ is bounded, we  have that 
$$
\text{$w_n\le u+1\in L^1(\Omega)$\quad for all $n\in \N$\quad and \quad$\lim_{n\to \infty}w_n=(u-1)^+$.}
$$ 
It follows  by the dominated convergent theorem we have
\begin{equation}\label{limt-Un-int}
\lim_{n\rightarrow\infty}U_n=  \int_{\Omega}|(u-1)^+(x)|^q\ dx.
\end{equation}
For $1<p<\infty$, note that  for every $v\in X_{s_1,p}\cap X_{s_2,q}$ and $x,y\in\R^N$ (see \cite{ML19}), it holds that
\begin{equation}\label{V+}
    |v(x)-v(y)|^{p-2}(v(x)-v(y))(v^+(x)-v^+(y))\ge |v^+(x)-v^+(y)|^p.
\end{equation}
  We can then  take $v=w_{n+1}\in X_{s_2,q}\cap X_{s_1,p}$ in \eqref{weakform} to get
\begin{align*}
&\cE_{s_1,p}(u,w_{n+1})+\cE_{s_2,q}(u,w_{n+1})=\lambda \int_{\Omega\cap\{w_{n+1}>0\}}|u|^{q-2}u~w_{n+1}\  dx.
\end{align*}
Therefore, it follows from \eqref{V+}  with  $v=u-C_n$, that
\begin{align*}
\|w_{n+1}\|^q_{{s_2,q}}\le &\iint_{\cQ}\frac{|w_{n+1}(x)-w_{n+1}(y)|^p}{|x-y|^{N+ps_1}}\  dxdy+\iint_{\cQ}\frac{|w_{n+1}(x)-w_{n+1}(y)|^p}{|x-y|^{N+ps_2}}  dxdy+U_{n+1} \\
&\le \cE_{s_1,p}(u,w_{n+1})+\cE_{s_2,q}(u,w_{n+1})+U_{n+1} \\
&=\lambda\int_{\Omega\cap\{w_{n+1}>0\}}|u|^{q-2}uw_{n+1}\ dx +U_{n+1}\\
&\le \big( \lambda(2^{n+1}-1)^{q-1}+1)U_n\le \lambda 2^{q(n+1)}U_n.
\end{align*}	
Moreover, we can observe that  $C_n=(1-{2^{-n}})=(1-\frac{2}{2^{n+1}})=C_{n+1}-\frac{1}{2^{n+1}}$, such that
\begin{align*}
w_{n+1}\le w_n=(u-C_n)^+= (u-C_{n+1}+\frac{1}{2^{n+1}})^+ = (w_{n+1}-[u-C_{n+1}]^-+\frac{1}{2^{n+1}})^+.
\end{align*}
As consequence, we have that~ 
$
\{w_{n+1}>0\}\subset \{w_{n+1}>\frac{1}{2^{n+1}}\}~
$
and 
\begin{equation}\label{measure}
|\Omega\cap\{w_{n+1}>0\}|\le  {2^{q(n+1)}}\int_{\Omega\cap\{w_{n+1}>\frac{1}{2^{n+1}}\}}|w_{n+1}(x)|^q\ dx\le {2^{q(n+1)}} U_n
\end{equation}
Now, using H\"older inequality and the fraction Sobolev inequality we have 
\begin{align*}
U_{n+1}=\int_{\Omega\cap\{w_{n+1}>0\}}|w_{n+1}(x)|^q\ dx
&\le |\Omega\cap\{w_{n+1}>0\}|^{\frac{qs_2}{N}}\|w_{n+1}\|^q_{L^{q^*_{s_2}
}(\Omega)}\\
&\le C|\Omega\cap\{w_{n+1}>0\}|^{\frac{qs_2}{N}}\|w_{n+1}\|_{s_2,q},
\end{align*}	
for some constant $C=C(N, s,\Omega)>0$ and all $n\in \N$. It thus follows from \eqref{measure} that
\begin{equation}\label{sequence-norm}
U_{n+1}\le\displaystyle C[2^{q(n+1)}]^{1+\frac{qs_2}{N}}U_n^{1+\frac{qs_2}{N}}=\tilde Cb^nU_n^{1+\alpha}\hspace{1cm}\text{for all $n\in \N$,}
\end{equation}
where $b:=2^{q(1+\frac{qs_2}{N})}>1$ and $\alpha:=\frac{qs_2}{N}$.  Note  in particular that  $w_0=(u-C_0)^+=u^+$ and 
\[
U_0=\int_{\Omega}|u^+|^q\ dx\le \int_{\Omega}|u|^q\ dx\le \delta^q.
\]
It follows from Lemma \ref{Giusti}  by conveniently taking $\delta$ such that $\delta<\tilde C^{-\frac{1}{2\alpha}}b^{-\frac{1}{2\alpha^2}}$ that 
\begin{equation}\label{limit-Un}
\lim_{n\rightarrow\infty}U_n = 0.
\end{equation}
Therefore, the limit in \eqref{limt-Un-int} yields that 
\[
\| (u-1)^+\|^q_{L^q(\Omega)}=0.
\]
Hence $(u-1)^+=0$ a.e. in $\Omega$, that is $u\leq 1$ a.e. in $\Omega$.  
By replacing $u$ by $-u$ in the above  computation, it follows that the rescaled $u$ satisfies
\begin{align*}
 \|u\|_{L^{\infty}(\Omega)}\le C\|u\|_{L^q(\Omega)},
\end{align*}
with some positive constant $C:=C(N,\Omega,s_1,s_2,q,p)$. We conclude that $u$ is bounded in $\Omega$. 

Next, for $x\in \R^N\setminus \Omega$, using the nonlocal $(p,q)$-Neumann condition we have that 
\begin{align*}
 |u(x)|&=   \left|\frac{C_{N,s}\displaystyle\int_{\Omega} \frac{
|u(x)-u(y)|^{r-2}u(y)}{|x-y|^{N+rs}}\ dy+C_{N,s}\int_{\Omega} \frac{
|u(x)-u(y)|^{r-2}u(y)}{|x-y|^{N+rs}}\ dy}{C_{N,s}\displaystyle\int_{\Omega} \frac{
|u(x)-u(y)|^{r-2}}{|x-y|^{N+rs}}\ dy+C_{N,s}\int_{\Omega} \frac{
|u(x)-u(y)|^{r-2}}{|x-y|^{N+rs}}\ dy}\right|\le \|u\|_{L^{\infty}(\Omega)}.
\end{align*}
Hence, $\|u\|_{L^\infty(\R^N\setminus\Omega)}\le \|u\|_{L^{\infty}(\Omega)}$. We conclude that $u$ is bounded in $\R^N,$  and the proof of Theorem \ref{Inf-bound} is completed. 
\end{proof}

\bibliographystyle{amsplain}

\begin{thebibliography}{10}

\bibitem{AM07} A. Ambrosetti and A. Malchiodi. Nonlinear \textit{analysis and semilinear elliptic problems}. Volume 104. Cambridge university press, (2007).

\bibitem{AFR22} A. Audrito,  J.C. Felipe-Navarro and X. Ros-Oton. \textit{The Neumann problem for the fractional Laplacian: regularity up to the boundary.} Annali scuola normale superiore-classe di scienze, (2022): 55-55.


\bibitem{BM19} L. Barbu and G. Moro{\c{s}}anu. \textit{Eigenvalues of the negative $(p,q)$-Laplacian under a Steklov-like boundary condition.} Complex Variables and Elliptic Equations 64, no. 4 (2019): 685-700.

\bibitem{BM2023} L. Barbu and G. Moro{\c{s}}anu. \textit{Full description of the spectrum of a Steklov-like eigenvalue problem involving the $(p,q)$-Laplacian.} Ann. Acad. Rom. Sci, Ser. Math. Appl (2023).

\bibitem{BM23} L. Barbu and G. Moro{\c{s}}anu. \textit{On the eigenvalue set of the $(p,q)$-Laplacian with a Neumann-Steklov boundary condition.} Differential and Integral Equations 36, no. 5/6 (2023): 437--452.

\bibitem{BM023} L. Barbu and G. Moro{\c{s}}anu. \textit{On eigenvalue problems governed by the $(p,q)$-Laplacian.} Studia Universitatis Babes-Bolyai, Mathematica 68, no. 1 (2023).

\bibitem{BM18} M. Bhakta and D. Mukherjee. \textit{Multiplicity results for $(p,q)$ fractional elliptic equations involving critical nonlinearities}. Advances in Differential Equations, 24(3/4), (2018): 185--228.

\bibitem{BDVV22} S. Biagi, S. Dipierro, E. Valdinoci, and E. Vecchi, \textit{Mixed local and nonlocal elliptic operators: regularity and maximum principles}, Communications in Partial Differential Equations 47, no. 3 (2022): 585--629.


\bibitem{BS22} N. Biswas and F. Sk. \textit{On generalized eigenvalue problems of fractional $(p, q)$-Laplace operator with two parameters.} Proceedings of the Royal Society of Edinburgh Section A: Mathematics (2022): 1--46.

\bibitem{BT24} V. Bobkov and M. Tanaka. \textit{Abstract multiplicity results for (p, q)-Laplace equations with two parameters.} Rendiconti del Circolo Matematico di Palermo Series 2 (2024): 1--28.

\bibitem{BLP14} L. Brasco, E. Lindgren and E. Parini. \textit{The fractional Cheeger problem.} Interfaces and Free Boundaries 16,  (2014): 419--458.


\bibitem{B11} H. Brezis and H. Br\'ezis. \textit{Functional analysis, Sobolev spaces and partial differential equations}. Vol. 2, no. 3. New York: Springer, (2011).



\bibitem{CEF24} C. Cowan, M. El Smaily, P. A. Feulefack, \textit{Existence and regularity results for a Neumann problem with mixed local and nonlocal diffusion.} Journal of Differential Equations, 423, (2025): 97--117.

\bibitem{DS15}  L. M. Del Pezzo and A. M. Salort. \text{The first non-zero Neumann p-fractional eigenvalue.} Nonlinear Analysis: Theory, Methods \& Applications 118 (2015): 130--143.

\bibitem{Valdinoci} E. Di Nezza, G. Palatucci and E. Valdinoci. \textit{Hitchhiker's guide to the fractional Sobolev spaces.} Bulletin des sciences math\'ematiques 136.5 (2012): 521--573.

\bibitem{SPV22} S. Dipierro, E. Proietti Lippi and E. Valdinoci, \textit{Linear theory for a mixed operator with Neumann conditions},  Asymptotic Analysis 128, no. 4 (2022): 571--594.

\bibitem{DV21} S. Dipierro and E. Valdinoci. \textit{Description of an ecological niche for a mixed local/nonlocal dispersal: an evolution equation and a new Neumann condition arising from the superposition of Brownian and L\'evy processes.} Physica A: Statistical Mechanics and its Applications 575 (2021): 126052.


\bibitem{DRV17} S. Dipierro,  X. Ros-Oton, and E. Valdinoci. \textit{Nonlocal problems with Neumann boundary conditions.} Revista Matem\'atica Iberoamericana 33, no. 2 (2017): 377--416.

\bibitem{FR78} E.R. Fadell, P.H. Rabinowitz, Generalized cohomological index theories for Lie group actions with an application to bifurcation questions for Hamiltonian systems. Inventiones mathematicae 45.2 (1978): 139--174.

\bibitem{FMS17} M. F{\u{a}}rc{\u{a}}{\c{s}}eanu, M. Mihailescu, and D. Stancu-Dumitru. \textit{Perturbed fractional eigenvalue problems.} Discrete and Contin. Dynam. Systems-A, 37(12): (2017): 6243--6255.

\bibitem{FMS15}M. F{\u{a}}rc{\u{a}}{\c{s}}eanu, M.  Mih{\u{a}}ilescu and D. Stancu-Dumitru. \textit{On the set of eigenvalues of some PDEs with homogeneous Neumann boundary condition.} Nonlinear Analysis: Theory, Methods \& Applications 116 (2015): 19--25.

\bibitem{Giusti03} E. Giusti. \emph{Direct methods in the calculus of variations}. World Scientific, (2003).


\bibitem{L23} G. Leoni. \textit{A first course in fractional Sobolev spaces}. Vol. 229. American Mathematical Society, (2023).

\bibitem{MM16} M. Mih{\u{a}}ilescu and G. Moro{\c{s}}anu. \textit{Eigenvalues of $-\Delta_p-\Delta_q$ under Neumann boundary condition.} Canadian Mathematical Bulletin 59, no. 3 (2016): 606--616.

\bibitem{MS15} P. Mironescu and W. Sickel. \textit{A Sobolev non embedding.} Atti Accad. Naz. Lincei Cl. Sci. Fis. Mat. Natur. Rend. Lincei (9) Mat. Appl. 26, no. 3 (2015): 291-298.

\bibitem{ML19} D. Mugnai and EP. Lippi. \textit{Neumann fractional p-Laplacian: Eigenvalues and existence results.} Nonlinear Analysis 188 (2019): 455--474.


\bibitem{SVWZ22} X. Su, E. Valdinoci, Y.  Wei, and J. Zhang, \textit{Regularity results for solutions of mixed local and nonlocal elliptic equations} Mathematische Zeitschrift 302, no. 3 (2022): 1855--1878.

\bibitem{ZF24} EWB. Zongo and P. A. Feulefack. \textit{Bifurcation results and multiple solutions for the fractional $(p, q)$-Laplace operators.}  Preprint arxiv:2406.15825 (2024).

\bibitem{ZR24} EWB. Zongo and B. Ruf. \textit{Bifurcation results for nonlinear eigenvalue problems involving the $(p,q)$-laplace operator}. To appear in Advances in Differential Equations (2024)

\end{thebibliography}

\end{document}